\documentclass[11pt,twoside,a4paper]{article}
\usepackage{amsmath,amssymb,amsthm}

\newtheorem{Thm}{Theorem}[section]
\newtheorem{Lem}[Thm]{Lemma}
\newtheorem{Pro}[Thm]{Proposition}
\newtheorem{Cor}[Thm]{Corollary}
\newtheorem{Def}[Thm]{Definition}

\theoremstyle{definition}

\theoremstyle{remark}
\newtheorem{Rem}[Thm]{Remark}

\newcommand{\Z}{\mathbb{Z}}
\newcommand{\N}{\mathbb{N}}
\renewcommand{\H}{\mathbb{H}}

\newcommand{\al}{\alpha}

\newcommand{\ga}{\gamma}
\newcommand{\Ga}{\Gamma}

\renewcommand{\phi}{\varphi}

\newcommand{\rank}{\operatorname{rank}}

\newcommand{\cone}{\operatorname{Cone}}

\newcommand{\ov}{\overline}

\newcommand{\dia}{\operatorname{dia}}

\begin{document}

\title{Embedding of hyperbolic Coxeter groups into products of binary trees and
aperiodic tilings}
\author{Alexander Dranishnikov\footnote{Supported by
NSF} \ \& Viktor Schroeder\footnote{Partially supported by SNF}}

\maketitle

\begin{abstract}
We prove that a finitely generated, right-angled, hyperbolic Coxeter group $\Ga$
can be quasiisometrically embedded into the product of
$n$ binary trees, where $n$ is the chromatic number of $\Ga$.
As application we obtain certain strongly aperiodic tilings
of the Davis complex of these groups.
\end{abstract}

\section{Introduction}

We consider a finitely generated, right-angled Coxeter group $\Gamma$ , 
i.e. a group
$\Gamma$ together with a finite set of 
generators $S$, such that every element of $S$ has order two and that all
relations in $\Gamma$ are consequences of relations of the form
$st = ts$, where $s,t \in S$.

We prove embedding results of the Cayley graph 
$C(\Gamma ,S)$ of hyperbolic right angled
Coxeter groups into products of the binary tree $T =T_{\{0,1\}}$.
On graphs and trees we consider always the simplicial metric, hence
every edge has length 1. On a product of trees we consider the
$l_1$-product metric, i.e. the distance is equal to the sum of
the distances in the factors.

\begin{Thm} \label{thm:mainthm}
Let $\Gamma$ be a finitely generated, hyperbolic, right-angled 
Coxeter group with chromatic number $n$. 
 Then there exists  a quasiisometric embedding 
$\psi : C(\Gamma,S) \to T\times \cdots \times T$ 
($n$-factors).

\end{Thm}

We can apply Theorem \ref{thm:mainthm} for special
Coxeter groups operating cocompactly on the hyperbolic plane $\H^2$, e.g.
the group generated by the reflections
at the edges of the regular right-angled geodesic hexagon.
This group has chromatic number $2$ and we obtain:

\begin{Cor} \label{cor:hypplane}
There exists a quasiisometric embedding of the hyperbolic plane into the product
of two binary trees.

\end{Cor}

Combining this with a result of Brady and Farb \cite{BF} we 
get the following higher
dimensional version:

\begin{Cor} \label{cor:higherdim}
For every integer $n$ there exists a quasiisometric embedding
$\psi: \H^n \rightarrow T\times\cdots\times T$ of the 
hyperbolic space $\H^n$ into the $2(n-1)$-fold product of binary trees.
\end{Cor}

It is an open problem, if for 
$n \geq 3$ there is a quasiisometric embedding of
$\H^n$ into the $n$-fold product of binary trees.
There are partial results in this direction. In
\cite{BS2} it is shown that there exists a quasiisometric embedding of
$\H^n$ into an $n$-fold product of locally infinite trees.
On the other hand
a recent construction of 
Januszkiewicz and 
Swiatkowski \cite{JS} shows for every $n$ the existence of
a right angled Gromov hyperbolic Coxeter group with virtual
cohomological dimension and coloring number equal to $n$.
Combining Theorem \ref{thm:mainthm} with that result
we obtain:

\begin{Cor} \label{cor:vcohdim}

For every
$n$ there exits a Gromov hyperbolic group $\Gamma_n$ with virtual
cohomological dimension $n$ which admits a quasiisometric embedding
into the product of $n$ binary trees.
\end{Cor}

Corollary \ref{cor:vcohdim} can be used to determine the hyperbolic rank
(compare \cite{BS1}) of a product of trees:

\begin{Cor} \label{cor:hyprank}
The hyperbolic rank of the product 
of
$n$ metric trees is $(n-1)$.
\end{Cor}

\begin{Cor} \label{cor:freegroups}
The hyperbolic rank of the product 
of
$n$ free groups is equal to $(n-1)$.
\end{Cor}

With the methods of the proof we also obtain
strongly aperiodic tilings of the Davis complex of right-angled,
hyperbolic Coxeter groups.

\begin{Thm} \label{thm:tiling-intro}
For every finitely generated, right-angled and hyperbolic
Coxeter group $\Ga$ the Davis complex $X$ admits a
strongly aperiodic tiling $\Phi$ with finitely many tiles.
In addition, every limit tiling of $\Phi$ is also strongly aperiodic.
\end{Thm}

In dimensions 2,3 and 4 there are right-angled reflection
groups operating with compact fundamental domain
on the hyperbolic space. Two dimensional examples come from the
regular right-angled $p$-gons, $p \geq 5$ in the hyperbolic plane.
In dimension 3 there exists a right-angled regular hyperbolic
dodecahedron (compare also \cite{A}). In dimension 4 there exists e.g.
a right-angled hyperbolic
120-cell ( \cite{PV},\cite{D2},\cite{C}).
Thus we obtain

\begin{Cor} \label{cor:tilinglowdim}
The hyperbolic spaces
$\H^2,\H^3$ and $\H^4$ admit strongly aperiodic tilings such that every limit
tiling is also strongly aperiodic
\end{Cor}

A strongly aperiodic tiling of $\H^2$ was recently constructed by
Goodman-Strauss \cite{G} (this tiling has even a finite strongly aperiodic set
of tiles). Existence of aperiodic tilings in
$\H^3$ and $\H^4$ follows from the results of 
Block-Weinberger \cite{BW}. The aperiodicity of their tilings
follow from the fact that the tilings are unbalanced.
We construct strictly balanced aperiodic tilings of the Davis complex.
Indeed we develop a general method how to modify tilings of the Davis
complex to obtain strictly balanced tilings.
For a precise statement see
Theorem \ref{thm:balanced}.

The structure of the paper is as follows.
After the preliminaries we prove the main result for the
two-colored case in section 3. In section 4 we generalize the argument to the
$n$-colored case. In section 5 we prove the corollaries of this result.
In section 6 we discuss aperiodic tilings of the Davis
complex.

It is a pleasure to thank Victor Bangert for a
useful discussion about the Morse Thue sequence.
Also we would like to thank the Max-Planck Institut 
f\"ur Mathematik in Bonn for the hospitality.

\section{Preliminaries}

\subsection{Right-angled Coxeter Groups} \label{subsec:racg}

In this section we review the necessary facts from the theory of
right-angled Coxeter groups.

A Coxeter matrix $(m_{s,t})_{s,t \in S}$ is a symmetric
$S \times S$ matrix with 1 on the diagonal and with all other
entries nonnegative integers different from 1.
A Coxeter matrix defines a Coxeter group $\Gamma$ generated by
the index set $S$ with relations
$(st)^{m_{s,t}} = 1$ for all $s,t \in S$.
Here we use the convention that $\gamma^0 =1$ for all elements, 
thus if $m_{s,t}=0$ then there
is no relation between $s$ and $t$. A Coxeter group $\Gamma$ 
is {\em finitely generated}, if
$S$ is finite. The group $\Gamma$ is called {\em right-angled}, if all
entries of the corresponding Coxeter matrix are 0,1,2. The 
Coxeter matrix of a right
angled Coxeter group is completely described by a graph 
with vertex set $S$ where we
connect two vertices $s$ and $t$ iff $m(s,t)=2$.

Let $\Gamma$ be a right-angled Coxeter group with generating set $S$.
The nerve $N=N(\Gamma, S)$ 
is the simplicial complex defined in the following way: the 
vertices of $N$ are the elements
of $S$. Two different vertices $s,t$ are joined by an edge, 
if and only if
$m(s,t)= 2$. In general (k+1) different vertices $s_1,\ldots ,s_{k+1}$ 
span a k-simplex, if and only
$m(s_i,s_k)=2$ for all pairs of different $i,j \in \{ 1,\ldots ,k+1\}$.

There is a simple characterization of hyperbolic right-angled
Coxeter groups.
A {\it square} in $S$ is a collection of
four elements
$t_1,t_2,s_1,s_2 \in S$, such that
$s_i$ commutes with
$t_1$ and $t_2$ for $i=1,2$, but
$s_1$ and $s_2$ as well as $t_1$ and $t_2$
do not commute. Thus
these elements are the vertices of a square in
the nerve.
The following characterization is known as the Siebenmann
no-square condition:

\begin{Lem} \label{lem:nosquare}
A right-angled Coxeter group $\Ga$ with generating set
$S$ is hyperbolic if and only if
$S$ contains no square.
\end{Lem}

The {\it chromatic number} of a 
graph is the minimal number of colors
needed to color the vertices 
in such a way that adjacent vertices have 
different colors. The chromatic number of a simplicial complex
is 
the chromatic number of its 1-dimensional skeleton.

Assume that 
the chromatic number of the nerve $N(\Gamma)$ of a right-angled
Coxeter 
group $\Gamma$ equals $n$ and let $c:N^{(0)}\to\{1,\dots, n\}$
be a 
corresponding coloring map. 
Then we can write $S$ as a disjoint union
$S = \bigsqcup_{c} S_{c}$ where 
$S_c$ are the elements of color $c$.
Elements of a given color do not commute.
We prefer to write
$S= A \bigsqcup B \bigsqcup C \ldots $, where $A \subset S$ are the
elements with color $a$.

Right-angled Coxeter groups have a very simple deletion law. By the following
two operations every word $W$ in the generators $S$ can be transformed 
to a reduced word and two reduced
words representing the same element can be transformed by means only the second
operation \cite{Br} :

\begin{description}
\item[(i)] delete a subword of the form $ss$, $s\in S$

\item[(ii)] replace a subword $st$ by $ts$ if $m_{s,t}=2$.
\end{description}

\noindent This deletion rule has the following consequences:

\begin{Lem} \label{deletionrule}

\begin{description}
\item[(a)] If $W$ and $W'$ represent the same element then the lenghts of
$W$ and $W'$ are either both even or both odd.

\item[(b)] Let $W$ and $W'$ be reduced representations of the same element
$\gamma \in \Gamma$, then $W$ and $W'$ are formed from the same set of
letters and they have the same length.

\end{description}

\end{Lem}

We now investigate some properties of the Cayley graph
$C(\Gamma,S)$ of $\Gamma$ with respect to the generating set $S$ for
a right-angled Coxeter group.
Let $\gamma \in \Gamma$ and let $W$ be a reduced representation of
$\gamma$. The length of $W$ is denoted by
$\ell (\gamma )$ and called the norm of $\gamma$. This is well defined by (b).
If $a$ is a color and
$\gamma \in \Ga$, then we denote with
$\ell_a(\ga)$ the number of letters with
color $a$ in a reduced representation of $\ga$.
This is well defined and we clearly have
$\sum_c \ell_c(\ga) = \ell(\ga)$.

On $\Gamma$ we consider the distance function
$d(\gamma,\beta) =\ell (\gamma^{-1} \beta )$.
Let $\gamma$ and $\beta$ be elements of $\Gamma$ which are
neighbors in the Cayley graph and let $W$ be a reduced word representing
$\gamma$. Then there exists a generator $s \in S$ such that
$\beta$ has the representation $Ws$. It follows from (b) that
$\ell (\gamma ) \neq \ell (\beta )$. Thus an edge in the Cayley graph
connects two elements with different norm.
As usual one can define geodesics in the Cayley graph.
A geodesic between two points $\alpha,\beta\in\Gamma$ is given by a
sequence
$\alpha=\gamma_0,\ldots,\gamma_k=\beta$
with $d(\gamma_i,\gamma_j)=\ \mid i-j\mid$.

Let $\alpha, \beta$ and $\gamma$ be elements of $\Gamma$. We say that
$\gamma$ {\it lies between} $\alpha$ {\it and} $\beta$ if
$d(\alpha,\gamma)  + d(\gamma,\beta) = d(\alpha,\beta)$.
The Cayley graphs of right angled Coxeter groups have the
following property, which says that any three
points in $\Gamma$ span a tripod:

\begin{Lem} \label{lem:between}
Let $\alpha,\beta,\gamma \in \Gamma$, then there exists $\delta \in \Gamma$
such that $\delta$ lies between $\gamma_i$ and $\gamma_j$ for
any choice of distinct elements $\gamma_i,\gamma_j \in \{ \alpha,\beta,\gamma\}$.
\end{Lem}

\begin{proof}

By the $\Gamma$-invariance of the metric it suffices to show this result
for the case that $\gamma = 1$.
Let $\alpha,\beta \in \Gamma$ and consider a geodesic path
$\alpha =\alpha_0,\ldots,\alpha_k=\beta$ from $\alpha$ to
$\beta$ and consider the sequence of norms
$n_0=\ell (\alpha_0 ),\ldots,n_k =\ell (\alpha_k )$. Note that 
by the properties discussed
above $\mid n_i - n_{i+1}\mid\ =1$.
Assume that there is a subsequence $\alpha_{i-1},\alpha_i,\alpha_{i+1}$ with
$n_{i-1}<n_i>n_{i+1}$.
Then one can represent $\alpha_i$ in a reduced way as $W_1t = W_2s$ where
$W_1$ is a reduced word representing $\alpha_{i-1}$ and
$W_2$ a reduced word representing $\alpha_{i+1}$ and $s,t \in S$.
Since by the deletion law one can transform $W_1t$ into $W_2s$ by means
of operations of type (b), we see that $st=ts$ and that one can represent
$\alpha_i$ by a reduced word of the form $Wst=Wts$ where $Ws$ represents
$\alpha_{i-1}$ and $Wt$ represents $\alpha_{i+1}$. Replace now $\alpha_i$ by
the element $\alpha_i'$
represented by $W$ and we obtain a new geodesic sequence between $\alpha$ and
$\beta$ such that for the corresponding sequence of norms
we have
$n_{i-1}>n_i'<n_{i+1}$.
Applying this procedure several times
we obtain a geodesic path from $\alpha$ to $\beta$, such that the sequence of the norms
has no local maximum any more and hence only a global minimum
$n_{i_0}$. The corresponding element $\delta = \alpha_{i_0}$ lies
between $\alpha$ and $\beta$ but also between
$1$ and $\alpha$ resp. $1$ and $\beta$.

\end{proof}

\subsection{Rooted Trees}

Let $\Omega$ be a set. We associate to $\Omega$ a rooted simplicial tree 
$T_{\Omega}$ in the following way
: the set of vertices is the set of finite sequences
$(\omega_1,\ldots ,\omega_k)$ with
$\omega_i \in \Omega$.
The empty sequence defines the root vertex and is denoted
by $\emptyset$.
Two vertices are connected by an edge in
$T_{\Omega}$ if their length (as sequences) differ by one
and the shorter can be obtained by erasing the last term of the longer.
The root vertex has $|\Omega|$ neighbors and every other vertex has
$|\Omega|+1$ neighbors, one ancestor and
$|\Omega|$ descendents. Here $|\Omega|$ denotes the
cardinality of $|\Omega|$.
The set $\Omega$ is also called the {\it label set} of
$T_{\Omega}$.
On trees we consider always the simplicial metric, hence
every edge has length 1. On a product of trees we consider the
$l_1$-product metric, i.e. the distance is equal to the sum of
the distances in the factors.

The tree $T_{\{0,1\}}$ is also called the 
{\it binary tree}.
If $\Omega$ is finite and $|\Omega| \ge 2$, then
$T_{\Omega}$ is quasiisometric to
the binary tree. This is easily proved using the fact that
the elements of $\Omega$ can be represented by binary sequences
of length $\le \log_2(|\Omega|)+1$.


\section{Two-colored Coxeter Groups} \label{sec:twocolored}

Since the case of chromatic number 2 is 
technically much easier, we first give the proof of
our main result in this case.
Thus we asssume in this section,
that $\Ga$ is generated by a finite set
$S$ which can be decomposed as
$S = A \bigsqcup B$ such that elements
within $A$ do not commute with each other and the same holds for
elements in $B$.

\subsection{Canonical $a$-presentation}

Fix a color, say $a$.
A 
{\it reduced left $a$-representation}
of $\ga \in \Ga$ is a reduced
word
$$W = W_1a_1W_2a_2\ldots W_ra_rW_{r+1}$$
representing
$\ga$
in which the
$W_i= w^i_{r_i}\ldots w^i_1$  are words with letters in $B$ and
the last entry of $W_i$ does not commute with
$a_i$, i.e. the letters with color $a$ are moved as left as possible. 
The reduced $a$-representation of
$\ga$ is unique and thus we call it also
the
{\it canonical $a$-presentation} of $\ga$.
The word
$W_i$
in the canonical $a$-representation is called the
{\it coefficient} of $a_i$ if
$i\le r$ and it is called the {\it free coefficient}
if
$i=r+1$.
Note that
$W_i$ could be empty.
Let $W$ be a reduced word, then the reduction to the
canonical $a$-representation permutes the letters
of $W$. We call this permutation the
{\it canonical reduction map} of $W$.

We now study in general the situation that
we have two 
canonical $a$-representations 
$U$ and $V$ such that the composition
$UV$ is a reduced word.
We choose the notation such that
$$U = U_1a_1 \cdots U_pa_pU_{p+1},$$ 
$$V= V_{p+1}a_{p+1} \cdots V_{p+m}a_{p+m}V_{p+m+1}.$$ 
Let then
$$W = W_1a_1W_2a_2\ldots W_{p+m}a_{p+m}W_{p+m+1}$$
be the canonical
$a$-presentation of the word $UV$ and let
$\phi:UV \to W$ be the canonical reduction map
for $UV$.
With this notation $\phi(a_i)=a_i$ for
all $1\le i\le p+m$.
Let
$U_{p+1}=u_q\ldots u_1$ be the free coefficient of
$U$ and $V_{p+1}=v_k\ldots v_1$ be the coefficient of
$a_{p+1}$ in $V$ and
let
$U_R=\{ u\in U | a_{p+1}< \phi(u)\}$ denote
the set of all letters $u \in U$ that
are moved by $\phi$ right of
$a_{p+1}$. Then

\begin{Pro} \label{pro:composition}
\begin{enumerate}
\item[(1)] $\phi(v_i)=w_i$ for all $i\le k$, where
$W_{p+1}=w_l\ldots w_1$ is the coefficient of $a_{p+1}$ 
in $W$.
\item[(2)] If $k\ge 1$ then $U_R = \emptyset$ and
$UV =W $ is the canonical $a$-representation.
\item[(3)] If $U_R \ne \emptyset$, then $U_R$ are the last $|U_R|$ letters of
$u_q\ldots u_1$ 
\item[(4)] If $\Ga$ is hyperbolic, then there is
at most one $u\in U_R$ with 
$\phi(u) > a_{p+2}$ . 
\end{enumerate}
\end{Pro} 

\begin{proof}
The first point is obvious.
If $k\geq 1$, then (since no element $u_i\in U_{p+1}$ can pass
$v_1$) $U_R = \emptyset$. This proves (2).
If $U_R \ne \emptyset$, then $k=0$ and clearly
$U_R$ consists out of the last letters of $U_{p+1}$, hence (3).
To prove the last point assume that there are
two elements in $U_R$ which are mapped by $\Phi$ to the right of
$a_{p+2}$. Then these two elements commute with
$a_{p+1}$ and $a_{p+2}$ and thus this four elements together form
a square. Hence $\Ga$ is not hyperbolic by Lemma \ref{lem:nosquare}.
\end{proof}

Consider two elements
$\ga, \ov{\ga} \in \Ga$.
We investigate the situation, how
the canonical $a$-representations of
these elements differ.
By Lemma \ref{lem:between}
$1,\ga$ and $\ov{\ga}$ span a tripod.
Thus there are words
$U$,$V$,$\ov{V}$ given in
canonical $a$-presentations, such that
$UV$ is a reduced representation of $\ga$ and
$U\ov{V}$ is a reduced representation of $\ov{\ga}$.
Furthermore
$V^{-1}\ov{V}$ is a reduced representation
of
$\ga^{-1}\ov{\ga}$.
Let $m=\ell_a(V)$ and
$\ov{m}=\ell_a(\ov{V})$.

We then write
$$U=U_1a_1\ldots U_pa_pU_{p+1}$$
$$V=V_{p+1}a_p \ldots V_{p+m}a_{p+m}V_{p+m+1}$$
$$\ov{V}=\ov{V}_{p+1}\ov{a}_p 
\ldots \ov{V}_{p+\ov{m}}\ov{a}_{p+\ov{m}}\ov{V}_{p+\ov{m}+1}$$

Let
$$W=W_1a_1\ldots  W_{p+m}a_{p+m}W_{p+m+1}$$
$$\ov{W}=\ov{W}_{1}\ov{a}_1 
\ldots \ov{W}_{p+\ov{m}}\ov{a}_{p+\ov{m}}\ov{W}_{p+\ov{m}+1}$$

be the canonical $a$-representations of
$\ga, \ov{\ga}$. Clearly
$W_ia_i=\ov{W}_i\ov{a}_i=U_ia_i$ for $1\le i\le p$.

\begin{Lem} \label{lem:pplusonedifferent}
$W_{p+1}a_{p+1} \ne \ov{W}_{p+1}\ov{a}_{p+1}$.
\end{Lem}

\begin{proof}

Since
$V^{-1}\ov{V}$ is reduced,
the first letter from
$V$ is different from the first letter
of $\ov{V}$.
We start by considering the case, that
the first letters are both of color
$a$, i.e.
$V_{p+1}=\emptyset=\ov{V}_{p+1}$.
In this case we have
$a_{p+1}\ne \ov{a}_{p+1}$ and hence the
claim.
If (say)
$V_{p+1}=\emptyset$ and
$\ov{V}_{p+1}\ne \emptyset$, then we get
from Proposition \ref{pro:composition}
that
$\ov{W}_{p+1}$ is longer then
$W_{p+1}$ and hence the claim.
If both
$V_{p+1} \ne \emptyset$ and 
$\ov{V}_{p+1} \ne \emptyset$,
then by Proposition \ref{pro:composition}
$W_{p+1}=U_{p+1}V_{p+1}$ and
$\ov{W}_{p+1}=\ov{U}_{p+1}\ov{V}_{p+1}$
and hence the two word differ at the first
entry of $V_{p+1}$ resp. $\ov{V}_{p+1}$. 
This proves the result.

\end{proof}

\subsection{Alice's diary}

Let $\kappa \in \N$.
Define
$E = (A \cup B \cup \{\emptyset\})$.
Let $T^{\dia}_a$ be the rooted tree
with label set $E^{\kappa}$.
We call $T^{\dia}_a$ the
{\it diary tree}.
We define a map
$\psi^{\dia}_a: \Ga \to T^{\dia}_a$ as follows:

Let $\ga \in \Ga$ be given by the canonical $a$-presentation
$$W=(w^1_{k_1}\ldots w^1_1)a_1\ldots (w^r_{k_r}\ldots w^r_1)a_r
(w^{r+1}_{k_{r+1}}\ldots w^{r+1}_1).$$

By $W|_i$ we denote the $i$-cut 
$$W|_i = (w^1_{k_1}\ldots w^1_1)a_1\ldots a_{i-1}(w^i_{k_i}\ldots w^i_1).$$

Recall that the
vertices of the tree
$T^{\dia}_a$  are finite sequences of elements
of $E^{\kappa}$.
We define
$\psi^{\dia}_a(W) = (\al_1,\ldots ,\al_i,\ldots ,\al_r)$
as a sequence of length $r$ by induction on $i$.

Let $\al_1$ be the string of the last
$\kappa$ symbols in the chain
$$\emptyset \ldots \emptyset  w^1_{k_1}\ldots w^1_1$$
considered as a word in the alphabet $E$.
We assume here that we have enough
(say $\kappa$) symbols $\emptyset$ in front.
We define
$\al_i$ as the string of the last $\kappa$ symbols
(E-letters) of the word in the alphabet $E$ that is obtained
from the word
$$\emptyset \ldots \emptyset  W|_i$$
by removing the E-letters from
$\cup_{j<i}\al_j$.

\begin{Rem}
One can intuitively describe this "diary" in the following way:
Consider the word $W$ as description of a long journey of
(say) Alice who
wants to write a diary about her trip. Every letter of $W$
represents a day.
There are two types of days during this journey.
The days of color $b$ where the weather is so fine
that Alice has no time to write her diary.
Then there are more cloudy days of color $a$. In the
morning of every $a$-day Alice writes $\kappa$ pages
in her diary: she starts with yesterday and then the day before
yesterday etc. Of course she skips the days which
were already described earlier in the diary. If there is no
day left to describe, she marks on the corresponding pages of her
diary the symbol $\emptyset$.
\end{Rem} 

It is not difficult to show and left as an
exercise to prove the following:
if
$\psi^{\dia}_a(W) = (\al_1,\ldots ,\al_i,\ldots ,\al_r)$
and
$\al_i$ contains the symbol $\emptyset$, then one
can reconstruct the whole i-cut
$W|_i$ from
$(\al_1,\ldots ,\al_i)$.

\begin{Rem} \label{rem:reconstruction}
It is a subtle point that Alice writes her diary in
the morning and not in the evening; i.e. that the entry
$\alpha_i$ does not contain $a_i$. This choice of the
diary makes it possible to reconstruct from the 
diary the i-cut
$W|_i$.
\end{Rem}

More generally we have

\begin{Pro} \label{pro:reconstruction}
Let
$W=W_1a_1\ldots W_{j+r}a_{j+r}W_{j+r+1}$ be a
canonical $a$-representation and assume that
$$k:= \kappa (r+1) - \ell(W_{j+1}a_{j+1}\ldots W_{j+r}) \geq 1.$$
Then one of the following holds:
\begin{enumerate}
\item[(1)] We can reconstruct the subword
$W_j a_j$ from the diary entries
$(\al_{j},\ldots ,\al_{j+r})$.
\item[(2)] The word $W_ja_j$ has $\ge k$ letters and
we can reconstruct from 
$(\al_{j},\ldots ,\al_{j+r})$
the last $k$ letters
of $W_j a_j$.
\end{enumerate}
\end{Pro}

\begin{proof}
If some entry
$\al_q$ for some $j\leq q\leq j+r$ contains
$\emptyset$, then we can reconstruct the complete $q$-cut
$W|_q$ from the diary entries
$(\al_1,\ldots ,\al_q)$.
In particular we can reconstruct in this case
$W_j a_j$ from the diary entries
$\al_s$, $j \leq s\leq j+r$.

Thus we can assume that the entries 
$\al_q$, $j\leq q\leq j+r$ 
do not contain
$\emptyset$ and hence each entry contains
$\kappa$ letters from the word
$W_1a_1\ldots a_{j+r-1}W_{j+r}$.
If we cannot reconstruct
$W_j a_j$ from these entries, then all
entries from
$\al_q$, $j\leq q\leq j+r$ are
from the word
$W_{j}a_{j}\ldots a_{q-1}W_q$
(if some entry is from $W_{j-1}a_{j-1}$, 
then the first occuring of
this entry is $a_{j-1}$ and in this moment we
have all information to reconstruct
$W_ja_j$).
Thus all
$\kappa (r+1)$ diary entries are from
$W_ja_j\ldots a_{j+r-1}W_{j+r}$, and hence at least
$k $
of them from
$W_j a_j$.

\end{proof}

We collect the maps
$\psi^{\dia}_a$ and the corresponding
$\psi^{\dia}_b$ to a common map
$$\psi^{\dia}:\Ga \to T^{\dia}_a \times T^{\dia}_b$$

On the product of the trees we consider the
$l_1$-product metric.

\begin{Lem} \label{lem:onelipschitz}
The map
$\psi^{\dia}$ is 1-Lipschitz.
\end{Lem}

\begin{proof}
Let $\ga, \ga' \in \Ga$ elements
with $d(\ga,\ga') =1$.
Then
$|\ell(\ga') -\ell(\ga)| =1$ and
we can assume that
$\ell(\ga')=\ell(\ga)+1$. In this
case we can represent
$\ga'$ in a reduced way as
$Us$, where $U$ is a reduced word
representing $\ga$. Assume w.l.o.g. that
$s\in A$. In this case clearly
$\psi^{\dia}_b(\ga)=\psi^{\dia}_b(\ga')$.
Note that
$\psi^{\dia}_a(\ga')$ has one more entry
as $\psi^{\dia}_a(\ga)$, and all other 
entries agree.
Thus 
$\psi^{\dia}(\ga)$ and $\psi^{\dia}(\ga')$ have
distance 1.

\end{proof}

The next Lemma shows that
the diary map is a "radial isometry".

\begin{Lem} \label{lem:radisometry}
The map $\psi^{\dia}$ is a radial isometry, i.e.
$|\psi^{\dia}(\ga)| = \ell(\ga)$,
where $|.|$ denotes the distance from
$(\mbox{root},\mbox{root}) \in T^{\dia}_a \times T^{\dia}_b$.
\end{Lem}

\begin{proof}
It follows from the definition of
$\psi^{\dia}_a$, that the sequence
$\psi^{\dia}_a(\ga)$ has $\ell_a(\ga)$ terms.
It follows that
$$|\psi^{\dia}(\ga)|= \ell_a(\ga)+ \ell_b(\ga)= \ell(\ga).$$
\end{proof}

A small disadvantage of the diary constructed above is
the fact that
$\al_i$ does not contain the information about
the last letter $a_i$ (compare remark \ref{rem:reconstruction}.
To obtain this additional information we define 
an "augmented" diary tree
$T'^{\dia}_a $ by the label set
$E^{\kappa} \times A$ and define

$$\psi'^{\dia}_a(\ga)=\psi'^{\dia}_a(W) = ((\al_1,a_1),\ldots ,(\al_r,a_r))$$
where $W$ is a canonical
$a$-representation of $\ga$.

The augmented diary map 
$\psi'^{\dia}=\psi'^{\dia}_a \times \psi'^{\dia}_b$ carries some
additional information.
Also this map is 1-Lipschitz and a radial isometry.

\begin{Rem} \label{rem:periodic}

The maps
$\psi^{\dia}$ and 
$ \psi'^{\dia}$ are (in most cases) not quasiisometric.
Consider reduced words
$\ov{A}$ formed from $a$-letters, and
$\ov{B}$ formed from $b$-letters, such that
the first and the last entry of $\ov{B}$ are different.
Then also the words $\ov{B}^k$ are reduced for $k \geq 1$.
Assume
$\ell(\ov{B}) \ll \ell(\ov{A})$.
Let $k \in \N$ such that 
$\ell(\ov{A}) \ll k$ and let also
$\kappa \ll k$.
Now consider the two reduced words
$\ga =(\ov{B})^k\ov{A}$ and
$\ov{\ga} =(\ov{B})^{k+1}\ov{A}$.
The conditions imply that
$\psi^{\dia}_a(\ga)=\psi^{\dia}_a(\ov{\ga})$ and
$$|\psi^{\dia}_b(\ov{\ga}) -\psi^{\dia}_b(\ga)| = \ell(\ov{B})
\ll 2\ell(\ov{A}) +\ell(\ov{B})= d(\ga,\ov{\ga}).$$
Thus $\psi^{\dia}$ cannot be quasiisometric.
It turns out that the essential point of this example
is the periodicity
within the words
$\ga$ and $\ov{\ga}$. To exclude this periodicity,
we use the Morse Thue decoration.
\end{Rem}

\subsection{The Morse Thue decoration}

The following sequence was studied by
Thue \cite{T} and later independently by Morse \cite{M}

\begin{Def}
(Morse-Thue Sequence t(i)).
Consider the substitution rule
$0\to 01$ and $1\to 10$. Then start from
$0$ to perform this substitutions
$$0\to 01 \to 0110\to 01101001\to \ldots $$
to obtain a nested family of sequences of length $2^k$ in the
alphabet $\{0,1\}$. 
The resulting limit
sequence is called the {\it Morse-Thue} sequence.
\end{Def}
 
The Morse-Thue sequence has the following remarkable
property (see e.g. \cite{He}).

\begin{Thm}
The Morse-Thue sequence contains no string of type
$WWW$ where $W$ is any word in $0$ and $1$.
\end{Thm}


\subsection{The quasiisometry}

We will use the Morse-Thue sequence to "decorate" elements
of the generating set.
To every generator $s\in S$ can be given a decoration
$0$ or $1$. Let
$E_a= \{b^0,b^1| b \in B\} \cup A \cup \{\emptyset\}$.
Thus $E_a$ is the alphabet $E$ where we in addition
have decorations of the letters of color $b$. 
Let
$T_a$ be the rooted tree with label set
$(E_a)^{\kappa} \times A$.

We define the final
map
$\psi_a$ from $\Ga$ into the tree $T_a$ in the following
way.
As above take the
canonical
$a$-representation $W$ of an element
$\ga \in \Ga$ and decorate every element
with color $b$ in $W$ by the Morse Thue sequence;
i.e. the $i$-th letter of color $b$ gets
the decoration $t(i)$. Then we apply the
map $\psi'^{\dia}_a(W)$ to the decorated word.

Analagously we define
$E_b$, $T_b$ and a map
$\psi_b:\Ga \to T_b$ by
interchanging the roles of the two colors.

Let
$$\psi= (\psi_a,\psi_b):\Ga \to T_a\times T_b$$

In the same way as above we see that
$\psi$ is 1-Lipschitz and a radial isometry.

\begin{Thm}
If $\Ga$ is  2-colored, right angled and
hyperbolic Coxeter group and let
$\kappa \ge 160$ then
the map $\psi$ is a quasiisometry.
\end{Thm}

\begin{proof}

To keep the proof more general, define
$n=2$ to be the chromatic number.

We know already that $\psi$ is 1-Lipschitz.
To prove the opposite inequality
let
$\ga, \ov{\ga} \in \Ga$
and let
$c:=d(\ga, \ov{\ga})$.

Claim: If
 $c \ge 60 n$ then
$|\psi(\ga)\psi(\ov{\ga})| \ge c/6n$.

Clearly the claim implies the desired lower estimate
on the distance of image points.
To prove the claim,
we use the notation as in
Lemma \ref{lem:pplusonedifferent}. 
Hence we have canonical
$a$-representations
$$U=U_1a_1\ldots U_pa_pU_{p+1}$$
$$V=V_{p+1}a_p \ldots V_{p+m}a_{p+m}V_{p+m+1}$$
$$\ov{V}=\ov{V}_{p+1}\ov{a}_p 
\ldots \ov{V}_{p+\ov{m}}\ov{a}_{p+\ov{m}}\ov{V}_{p+\ov{m}+1}$$
such that
$$W=W_1a_1\ldots  W_{p+m}a_{p+m}W_{p+m+1}$$
represents $\ga$ and $\ov{\ga}$ is represented by
$$\ov{W}=\ov{W}_{1}\ov{a}_1 
\ldots \ov{W}_{p+\ov{m}}\ov{a}_{p+\ov{m}}\ov{W}_{p+\ov{m}+1}$$

We can assume without loss of
generality that in a
reduced representation of 
$\ga^{-1}\ov{\ga}$ the number of 
letters with color $a$ is larger or
equal than the number of letters with
color $b$.

This means that
$m+\ov{m} \ge c/n$.

We first show that the claim is true in the
case 
$m < c/3n$ or
$\ov{m} < c/3n$:
Assume that (say)
$m < c/3n$.
Since
$m + \ov{m} \ge c/n$
we have
$\ov{m}-m\ge c/3n$.
Since
$\psi$ is 
a radial isometry, we have

$$|\psi_a(\ov{\ga})|-|\psi_a(\ga) |= 
(p+\ov{m})-(p+m)=\ov{m}-m \ge c/3n$$
and hence
$|\psi_a(\ga)\psi_a(\ov{\ga})|\ge c/3n$.
Thus we can assume for the rest of the proof that
$m \ge c/3n$
and
$\ov{m} \ge c/3n$.

We will now prove the claim by contradiction.
Therefore assume that
$|\psi_a(\ga)\psi_a(\ov{\ga})| < c/6n$.
This assuption together with
$m \ge c/3n$
and
$\ov{m} \ge c/3n$
implies

\begin{equation} \label{eq:psi}
\psi_a(W|_{p+r} a_{p+r} )
=\psi_a(\ov{W}|_{p+r} \ov{a}_{p+r}) \ \ \mbox{for some}\ \ r> c/6n
\end{equation}

In particular the diary entries
$\psi^{\dia}_a(W|_{p+r} a_{p+r})=(\al_1,\ldots ,\al_{p+r})$ coincide
with the corresponding entries
$(\ov{\al}_1,\ldots ,\ov{\al}_{p+r})$ .
Since
$c \ge 60 n$ by assumption we have
$r > 10$ and one computes elementary that also
$(r-1)\ge c/10n$ and $(r-2)\ge c/10n$. 

Consider now the three cases :

\smallskip

\noindent Case 1: $V_{p+1}=\emptyset=\ov{V}_{p+1}$.
Then it follows from the fact that
$V^{-1}\ov{V}$ is reduced that
$a_{p+1}\ne \ov{a}_{p+1}$
and hence
$\psi_a(W|_{p+1} a_{p+1} )
\ne \psi_a(\ov{W}|_{p+1} \ov{a}_{p+1})$
(by definition of the augmented map $\psi'^{\dia}_a$).
This is a contradiction to 
equation (\ref{eq:psi}).

\smallskip

\noindent Case 2: $V_{p+1}=\emptyset$ and
$\ov{V}_{p+1}\ne \emptyset$.

In this case
$\ov{W}_j = \ov{V}_j$ for all $j\ge p+2$ by
Proposition \ref{pro:composition} (2).
In particular
$$\ell(\ov{W}_{p+3}\ov{a}_{p+3}\ldots \ov{W}_{p+r}) \le 
\ell(\ov{V}) \le c$$
and hence (since $\kappa \ge 20 n$ and $(r-2)\ge c/10n$).
$$k = \kappa (r-2) -\ell(\ov{W}_{p+3}\ov{a}_{p+3}\ldots \ov{W}_{p+r})
 \ge   c $$ 
By Proposition \ref{pro:reconstruction}
we can reconstruct from
$(\ov{\al}_{p+2},\ldots ,\ov{\al}_{p+r})$ either the
whole word
$\ov{W}_{p+2}\ov{a}_{p+2}$ or at least the last
$k\ge c$ letters of this word.
Since
$\ov{W}_{p+2}\ov{a}_{p+2}$ is a subword
of $\ov{V}$ it has length $\le c$ and hence
we can reconstruct the whole word.
Then, since
$(\ov{\al}_{p+2},\ldots ,\ov{\al}_{p+r})=
(\al_{p+2},\ldots ,\al_{p+r})$ we conclude
$\ov{W}_{p+2}\ov{a}_{p+2} = W_{p+2}a_{p+2}$.

We make a similar computation 
at the place $p+1$.
In the same way as above we compute 
$$\ell(\ov{W}_{p+2}\ov{a}_{p+2}\ldots \ov{W}_{p+r}) \le \ell(\ov{V}) \le c$$
and hence (since $\kappa \ge 80 n$ and $(r-1) \ge c/10n$)
$$k = \kappa (r-1) -\ell(\ov{W}_{p+2}\ov{a}_{p+2}\ldots \ov{W}_{p+r})
 \ge  7c .$$ 
By Proposition \ref{pro:reconstruction}
we can reconstruct from
$(\ov{\al}_{p+1},\ldots ,\ov{\al}_{p+r})$ either the
whole word
$\ov{W}_{p+1}\ov{a}_{p+1}$ or at least the last
$k$ letters of this word.
Assume for a moment that we can reconstruct the whole word.
Then, since
$(\ov{\al}_{p+1},\ldots ,\ov{\al}_{p+r})=
(\al_{p+1},\ldots ,\al_{p+r})$ we have
$\ov{W}_{p+1}\ov{a}_{p+1} = W_{p+1}a_{p+1}$ in
contradiction to Lemma \ref{lem:pplusonedifferent}.

It follows that
$\ov{W}_{p+1}\ov{a}_{p+1}$ and
$W_{p+1}a_{p+1}$ have at least
$k\ge 7c$ letters and the last
$k$ letters of the words coincide.

Note that (using the terminology of Proposition
\ref{pro:composition}):
$\ov{W}_{p+1} = W_{p+1}U_R\ov{V}_{p+1}$.
By the last point of Proposition \ref{pro:composition}
and the just proved fact that $\ov{W}_{p+2}=W_{p+2}$ we see
$$|U_R| \le |W_{p+2}|+1 =|\ov{W}_{p+2}|+1 \le \ell(\ov{V}) \le c.$$
It follows that
$|U_R\ov{V}_{p+1}| \le 2c$.
Since the last $6c$ letters of
$\ov{W}_{p+1}$ and $W_{p+1}$ coincide, 
we can write
$W_{p+2}$ in the form  $MHH$ and
$\ov{W}_{p+2}$ in the form $MHHH$
where
$H=U_R\ov{V}_{p+1}$.

Since $\psi_a$ also contains the Morse-Thue
decoration of the elements of $B$ and
the property of the Morse-Thue sequence, there can not be
any subsequence of decorated letters of the form $HHH$ . 
Thus
we obtain a contradiction.

\smallskip

\noindent Case 3. $V_{p+1} \ne \emptyset$ and 
$\ov{V}_{p+1} \ne \emptyset$.

Then
$\ov{W}_{p+1}=U_{p+1}\ov{V}_{p+1}$ and
$W_{p+1}=U_{p+1}V_{p+1}$.
Since 
$\ov{V}_{p+1} \ne V_{p+1}$ but
the last $7c$ letters of
$\ov{W}_{p+1}$ and $W_{p+1}$ coincide by the same
arguments as in Case 2,
we see
$\ell(\ov{V}_{p+1}) \ne \ell(V_{p+1})$.
Let w.l.o.g.
$\ell(\ov{V}_{p+1}) > \ell(V_{p+1})$.
It is now completely elementary to see that again
$W_{p+2}$ is of the form  $MHHP$ and
$\ov{W}_{p+2}$ of the form $MHHHP$
where $H$ are the first
$\ell(\ov{V}_{p+1}) - \ell(V_{p+1})$ letters of $\ov{V}_{p+1}$.
We obtain now a contradiction as in Case 2.

\end{proof}

\section{The n-colored case}

In this section we modify the arguments of
section \ref{sec:twocolored} to prove the general case.
We assume that
$\Ga$ is a finitely generated right-angled Coxeter group
with chromatic number $n$.
Thus the set $S$
of generators has a decomposition
$S = A \bigsqcup B \bigsqcup \ldots$ into $n$ subsets in such a way that
elements of the same color do not commute.
Fix a color, say $a$. 
A {\it reduced left $a$-representation} of
$\ga \in \Ga$ is a recudes word $W$ representing
$\ga$, which has the form
$$W=W_1a_1W_1a_2\ldots W_ra_rW_{r+1}$$
in which for every $i$ no letter  $x$ in $W_i$
commutes with $a_i$ and with all letters of
$W_i$ to the right of $x$.

As discussed in
\ref{subsec:racg} every set of mutually commutative generators
$R \subset S$ defines a simplex in the
nerve $N(\Ga,S)$. For every element
$\ga\in \Ga$ we define a 
{\it right simplex presentation}
$$U=\Delta_q\ldots \Delta_1$$
where the letters of each word $\Delta_i$ form a simplex.
The simplices $\Delta_i$ are defined by induction
on $i$. We consider $U$ in reduced presentation
and define $\Delta_1$ as the set of all letters
from $U$ which can be placed at the very end of the word.
So $\Delta_1$ consists of all letters $x$ that 
commute with all letters to the right.
Then $U$ can be represented as $U^1\Delta_1$. Then we apply
this procedure to $U^1$ to obtain
$\Delta_2$ and so on.
We note that the right simplex presentation
is unique in the sense that the sequence of simplices
is uniquely determined.
Thus, it gives a unique presentation of a group
element up to permutation
of the letters inside simplices.
A right simplex presentation $U=\Delta_q\ldots \Delta_1$ is
a word in the alphabet $S$. On the other hand it
is a word in the alphabet 
$\Sigma$, the set of all simplices. We will
use the same notation for both.

We now define a {\it canonical $a$-representation} of
$\ga \in \Ga$ to be a reduced left $a$-presentation
$$W=W_1a_1W_1a_2\ldots W_ra_rW_{r+1}$$
where each
$W_i=\Delta^i_{r_i}\ldots \Delta^i_1$ is in the right
simplex presentation.

The word
$W_i$
in the left $a$-representation is called the
{\it coefficient} of $a_i$ if
$i\le r$ and it is called the free coefficient,
if
$i=r+1$.

Let
$\Sigma_a$ be the set of nontrivial
simplices of $N(\Ga,S)$ which do not contain
letters with color $a$.
We can view $W$ as a word in the
alphabet
$\cal{A}$  $= A \cup \Sigma_a$. Considered as a word
in this alphabet it is unique and therefore also
called
{\it the} canonical $a$-representation.
Note that considered as a word in $S$ a simplex from
$\Sigma_a$ has length $\le n-1$ since it cannot contain
letters with the same color.

Let $W$ be a reduced word, a reduction to a
canonical $a$-representation permutes the letters
of $W$. We call this permutation the
{\it canonical reduction map} of $W$.

As above we study the situation that
we have 
$$U = U_1a_1 \cdots U_pa_pU_{p+1},$$ 
$$V= V_{p+1}a_{p+1} \cdots V_{p+m}a_{p+m}V_{p+m+1}.$$ 
and
$$W = W_1a_1W_2a_2\ldots W_{p+m}a_{p+m}W_{p+m+1}$$
the canonical
$a$-presentation of the word $UV$ and let
$\phi:UV \to W$ be the canonical reduction map
for $UV$.
With this notation $\phi(a_i)=a_i$ for
all $1\le i\le p+m$.
Let
$U_{p+1}=\Delta^u_q\ldots \Delta^u_1$ be the free coefficient of
$U$ and $V_{p+1}=\Delta^v_k\ldots \Delta^v_1$ be the coefficient of
$a_{p+1}$ in $V$ and
let
$U_R=\{ u\in U | a_{p+1}< \phi(u)\}$ denote
the set of all letters $u \in U$ that
are moved by $\phi$ right of
$a_{p+1}$. Then

\begin{Pro} \label{pro:composition:n}
\begin{enumerate}
\item[(1)] $\phi(\Delta^v_k) \subset \phi(\Delta^w_k)$ where
$\Delta^w_k\ldots \Delta^w_1= W_{p+1}$ is the coefficient
of $a_{p+1}$ of $W$.
\item[(2)] If $\Ga$ is hyperbolic, then the following holds:
If $u \in U_R$ and $\phi(u) > a_{p+2}$, then $u \in \Delta^u_1$. 
\item[(3)]  If $\Ga$ is hyperbolic, and $V_{p+1}\ne \emptyset$, then
$U_R \subset \Delta^u_1$. 

\end{enumerate}
\end{Pro}

\begin{proof}
(1). Let
$U_{p+1}$ be the free coefficient of
$U$ and $V_{p+1}$ be the first coefficient of
$V$. By the definition of the canonical
$a$-reduction the permutation
$\Phi$ is a composition
$\Phi=\Phi_2\circ\Phi_1$ of two
maps, where $\Phi_1$ moves
the $a$'s to the left as far as possible.
In other words
$\Phi_1$ moves $U_R$ to the right of
$a_{p+1}$.
Then $\Phi_2$ forms the simplices
$\Delta^w_k\ldots \Delta^w_1= W_{p+1}$ out
of $\ov{U}_{p+1}V_{p+1}$ where
$\ov{U}_{p+1}=U_{p+1} \setminus U_R$.
One now checks easily
$\Phi(\Delta^v_k) \subset \Phi(\Delta^w_k)$.

(2). Assume that there exists
$u \in U_{p+1} \setminus \Delta^u_1$ with
$\phi(u) > a_{p+2}$. Since $u$ cannot go right
of all elements of $\Delta^u_1$, there exists
$u' \in \Delta^u_1$ which does not commute with
$u$. Then $\phi(u')>\phi(u) >a_{p+2}$.
If $a_{p+1}\ne a_{p+2}$ then
$u,u',a_{p+1},a_{p+2}$ form a square in contradiction
to the hyperbolicity assumption.
If $a_{p+1}= a_{p+2}$, then there exists $u''$ between 
$a_{p+1}$ and $a_{p+2}$ which does not commute with
them. Then $u,u',u'',a_{p+1}$ form a square.

(3). If we assume the contrary, then similar to (2)
we obtain a square formed by two elements from
$U_R$, one element from $V_{p+1}$ and $a_{p+1}$.
\end{proof}

In the same way as in section \ref{sec:twocolored}
represent two given elements
$\ga, \ov{\ga} \in \Ga$ by means of words
$$U=U_1a_1\ldots U_pa_pU_{p+1}$$
$$V=V_{p+1}a_p \ldots V_{p+m}a_{p+m}V_{p+m+1}$$
$$\ov{V}=\ov{V}_{p+1}\ov{a}_p 
\ldots \ov{V}_{p+\ov{m}}\ov{a}_{p+\ov{m}}\ov{V}_{p+\ov{m}+1}$$
 such that $\ga$ and $\ov{\ga}$ are represented by
$$W=W_1a_1\ldots  W_{p+m}a_{p+m}W_{p+m+1}$$
$$\ov{W}=\ov{W}_{1}\ov{a}_1 
\ldots \ov{W}_{p+\ov{m}}\ov{a}_{p+\ov{m}}\ov{W}_{p+\ov{m}+1}$$

Clearly
$W_ia_i=\ov{W}_i\ov{a}_i=U_ia_i$ for $1\le i\le p$ and we have
in generalization of Lemma \ref{lem:pplusonedifferent}

\begin{Lem} \label{lem:pplusonedifferent:n}
$W_{p+1}a_{p+1} \ne \ov{W}_{p+1}\ov{a}_{p+1}$.
\end{Lem}

\begin{proof}

Since
$V^{-1}\ov{V}$ is reduced,
the first letter from
$V$ is different from the first letter
of $\ov{V}$.

\smallskip

\noindent Case 1:$V_{p+1}=\emptyset=\ov{V}_{p+1}$.
In this case we have
$a_{p+1}\ne \ov{a}_{p+1}$ and hence the
claim.

\smallskip

\noindent Case 2: $V_{p+1}=\emptyset$ and
$\ov{V}_{p+1}\ne \emptyset$.

We can assume in addition that
$a_{p+1} =\ov{a}_{p+1}$ since otherwise the result is
trivially true. Note
that no element from
$\ov{V}_{p+1}$ can be moved right to $\ov{a}_{p+1} = a_{p+1}$.
Therefore all elements in $U_{p+1}\ov{V}_{p+1}$ which can
be moved right of $\ov{a}_{p+1}$ are elements of
$U_{p+1}$ and hence elements from
$U_R$ in the notation of Proposition \ref{pro:composition:n}
If follows that $\ell(W_{p+1}) <\ell(\ov{W}_{p+1})$.

\smallskip

\noindent Case 3: $V_{p+1} \ne \emptyset$ and 
$\ov{V}_{p+1} \ne \emptyset$.
 As above we can assume in addition that
$a_{p+1} = \ov{a}_{p+1}$.

We first investigate the case
$U_R \ne U_{\ov{R}}$.
Assume that there exists $u \in U_R$
such that
$u$ is not in $U_{\ov{R}}$.
Then there exists a "blocking element"
$b \in \ov{V}_{p+1}$ which does not
commute with $u$.
The letter $b$ is not contained in
$V_{p+1}$, since $u \in U_R$.
Furthermore $b$ is not contained
in $U_{\ov{R}}$ since
$U_{\ov{R}}\cap \ov{V}_{p+1} = \emptyset$.
It follows that the letter $b$ occurs
more often in
$\ov{W}_{p+1}$ than in $W_{p+1}$. But this
implies the claim.

Thus we can assume that
$U_R = U_{\ov{R}}$.
Now by Lemma \ref{pro:composition:n}(3)
$U_R =U_{\ov{R}} \subset \Delta^u_1$ which implies
that we can write $U_{p+1}$ in a reduced representation as
$U_{p+1} = U'U_R = U'U_{\ov{R}}$.
Hence the word
$W_{p+1}$ represents the same group element as
$U'V_{p+1}$ and
$\ov{W}_{p+1}$ the same element as
$U'\ov{V}_{p+1}$.
Since 
$V_{p+1}$ represents a different group element as
$\ov{V}_{p+1}$ we obtain the result.

\end{proof}

We can now define the diary map.
Let $\kappa \in \N$.
Define
$E_a = ( \cal{A} \cup \{\emptyset\})$.
Then
$E_a^{\kappa}$ is the label set $T^{\dia}_a$, the
{\it diary tree}.

In the same way as in section \ref{sec:twocolored} we define the
diary map

$$\psi^{\dia}_a(W) = (\al_1,\ldots ,\al_i,\ldots ,\al_r)$$
as a sequence of length $\ell_a(W)$. 
As above we obtain 

\begin{Pro} \label{pro:reconstruction:n}
Let
$W=W_1a_1\ldots W_{i+m}a_{i+m}W_{i+m+1}$ be a
canonical $a$-representation and assume that
$$k:= \kappa m - \ell_{\cal{A}}(W_{i+1}a_{i+1}\ldots W_{i+m}) \geq 1,$$
where $\ell_{\cal{A}}$ denotes the length of a word in the alphabet
$\cal{A}$.
Then one of the following holds:
\begin{enumerate}
\item we can reconstruct
$W_i a_i$ from the diary entries
$(\al_{i+1},\ldots ,\al_{i+m})$
\item The word $W_i a_i$ has $\ge k$ $\cal{A}$-letters and
we can reconstruct from 
$(\al_{i+1},\ldots ,\al_{i+m})$
the last $k$ $\cal{A}$-letters
of $W_i a_i$.
\end{enumerate}
\end{Pro}

We collect the maps
$\psi^{\dia}_c$ to a common map
$\psi^{\dia}:\Ga \to \prod_c T^{\dia}_c$.
and we also define the augmented diary map.
As in chapter \ref{sec:twocolored} we see that the diary map is
Lipschitz and a radial isometry.
Finally we decorate the simplices from
$\Sigma_a$ by the Morse Thue sequence in the following way.
Given a word $W$ in the canonical $a$-representation,
we decorate for every color
$b \ne a$ every $S$-letter with color $b$ by
the Morse-Thue sequence, i.e. the $i$-th letter of
color $b$ gets the decoration $t(i)$.
Since within a given color the order of the elements 
cannot change, this is well defined.
Thus a simplex is now a simplex of decorated $S$-letters.
To this decorated sequence we apply the augmented
diary map to obtain the map
$\psi_a$ and finally we collect the maps to get
$$\psi= \prod_c \psi_c :\Ga \to \prod_c T_c$$

Now we can prove
Theorem \ref{thm:mainthm}
in general. The proof is essentially
as in section \ref{sec:twocolored}.
There are however some significant differences in the discussion
of Case 2 and Case 3.

The idea is to reduce the argument to the case of
two colors.
Therefore we use the following reduction.
Let $W$ be a word given in the canonical $a$-decomposition and
let $b$ be a color different from $a$.
We define a word
$W^b$ in the alphabet
$\{b^o, b^1 | b\in B\} \cup A \cup \{\star\}$ 
(where $\star$ is some additional symbol) in the 
following way:
Replace every decorated simplex $\Delta \in \Sigma_a$
from $W$ by the decorated letter with color $b$, if 
$\Delta$ contains a letter with color $b$.
If $\Delta$ does not contain a letter with color
$b$ then replace $\Delta$ by $\star$. 

\begin{Lem} \label{lem:star}
Let
$G$ and $H$ be words given
in the canonical $a$-representation and assume
that $GH$ is reduced.
Assume that the first $\mathcal{A}$ letter of
$H$ contains a letter of color $b$.
Then $(GH)^b$ can be obtained
from $G^b H^b$ by removing certain
numbers of the letter $\star$ from
$G^b$.
\end{Lem}

\begin{proof}
Let $H = \Delta_l\ldots \Delta_1$ be the
canonical $a$-decomposition of $H$.
Let $b$ be a color that occurs in $\Delta_l$.
Since $GH$ is reduced, there is no letter
of color $b$ in $G$ which is mapped by the 
canonical reduction map right to a letter of color
$b$ in $\Delta_l$. The result follows now easily.
\end{proof}

We now proceed with the proof of the main theorem in the 
$n$-colored case.

Case 2: $V_{p+1} =\emptyset$ and $\ov{V}_{p+1} \ne \emptyset$.

Similar in chapter \ref{sec:twocolored} we conclude
that
$\ov{W}_{p+1} a_{p+1}$ and 
$W_{p+1}a_{p+1}$ have the last
$k\ge 7c$ letters in common. These
letters are decorated simplices.
We furthermore see that in this case we can write
$\ov{W}_{p+1} = W_{p+1}H$ where $H = U'\ov{V}_{p+1}$
and $U'=U_R \setminus \ov{U}_R$.
Let $H = \Delta_l\ldots \Delta_1$ be the
canonical $a$-decomposition of $H$.
Let $b$ be a color that occurs in $\Delta_l$.
This implies that there is no letter with color $b$
in $\ov{U}_R$.
Let $\Delta^w_l\ldots \Delta^w_1$ be the last
simplices of $W_{p+1}$. Since the last $\cal{A}$-letters of
$W_{p+1}$ and $\ov{W}_{p+1} = W_{p+1}H$ 
coincide we have
\begin{equation} \label{eq:subs}
\Delta_i \subset \Delta^w_i 
\end{equation}
by
Proposition \ref{pro:composition:n} (1).

We use now the reduction to the 2-colored case.
By Lemma \ref{lem:star}
we obtain
$(W_{p+1}H)^b$ from
$W^b_{p+1} H^b$ by deleting some $\star$-letters from
$W^b_{p+1}$.

Because of relation (\ref{eq:subs}) and
since $\ov{U}_R$ contains no letter of color
$b$, we see
that the last $k \ge 7c$ letters of
$W^b_{p+1}$ coincide with 
$H^b$. 
It follows from Lemma \ref{lem:star} that
there are words $Q, H_1, H_2 \in \cal{A}$, such that
$W_{p+1}=QH_1$,
$\ov{W}_{p+1}= QH_1H_0$ such that
$H_0^b=H^b$ and
$H_1^b$ can be obtained from
$H^b$ by deleting some $\star$-letters.
In particular the first letter of
$H_1^b$ has color $b$.
Again, since the
last letters of
$W^b_{p+1}=(QH_1)^b$ coincide with the last letters from
$\ov{W}^b_{p+1} = (QH_1H_0)^b$ we see that
$\ov{W}^b_{p+1} = \sigma H^b_2H^b_1H^b$ where
$H^b_2$ is obtained from
$H^b_1$ by deleting some $\star$-letters.
Thus the sequence of simplices with a vertex of color
$b$ has a subsequence of period 3 in contradiction to
the properties of the decoration.

Case 3: $V_{p+1} \ne \emptyset$ and $\ov{V}_{p+1} \ne \emptyset$.

We first show, that the length of 
$V_{p+1}$ differs from the length of
$\ov{V}_{p+1}$ where the length is measured in
the alphabet $\cal{A}$.

Assume to the contrary that
$$V_{p+1}=\Delta^v_l \ldots \Delta^v_1$$
and
$$\ov{V}_{p+1}=\Delta^{\ov{v}}_l \ldots \Delta^{\ov{v}}_1$$
are of the same length.

Since
$V\ov{V}^{-1}$ is reduced,
$\Delta^{\ov{v}}_l \cap \Delta^v_l = \emptyset$.
Let 
$x \in \Delta^{\ov{v}}_l$ and
$y \in \Delta^v_l$ letters (of $S$) . 

Since the last $7c$ $\cal{A}$-letters of
$W_{p+1}$ and $\ov{W}_{p+1}$ coincide, we see
$\Delta^w_l =\Delta^{\ov{w}}_l$ 
and by
Proposition \ref{pro:composition:n} (1)
$\Delta^{\ov{v}}_l \subset\Delta^w_l$ and
$\Delta^{v}_l \subset\Delta^w_l$. Thus the
letters $x,y$ must also occur
in $U$ in a way, that the letter
$x$ can be moved in $W$ to
$\Delta^w_l$ and $y$ in $\ov{W}$
to $\Delta^{\ov{w}}_l$.
Since $U\ov{V}$ is reduced, there exists
$z \in U$ right to $x$ which does not commute with
$x$. Since this $x$ can be moved to $\Delta^v_l$ we see
that also $z$ can be moved right to $\Delta^v_l$ which implies
that $y$ commutes with $z$.
In the same way we see that there exists $u \in U$ right of
$y$ which does not commute with $y$ but commutes with $x$.
The letters $x,y,z,u$ form a square in contradiction to
the hyperbolicity condition.

We therefore can assume w.l.o.g. that
$\ov{V}_{p+1}$ is longer than $V_{p+1}$.
Let $b$ be a color which accurs in 
the first $\cal{A}$-letter of 
$\ov{V}_{p+1}$.
With arguments as in case 2 one shows now that
$\ov{W}^b_{p+1}$ is obtained from a word of the form
$M^bH^b_2H^b_1H_0^bP^b$ by removing some $\star$-letters.
Here $H^b_i$ is obtained from $H^b$ by removing
$\star$-letters and 
$H$ is formed by the first 
$\ell_{\cal{A}}(\ov{V}_{p+1}) -\ell_{\cal{A}}(V_{p+1})$
$\cal{A}$-letters of $\ov{V}_{p+1}$.
We obtain a contradiction in a similar way.

\section{Proof of the Corollaries }

\begin{proof} (of Corollary \ref{cor:hypplane})

Consider the right angled Coxeter group
$\Gamma$ given by the generator set $S=\{s_1,\ldots,s_6\}$
and relations $s_i s_{i+1} = s_{i+1} s_i $ (indices mod 6).
This group acts discretely on the hyperbolic plane 
such that a Dirichlet fundamental
domain is bounded by the regular right angled hexagon in $\H^2$.
By
Theorem \ref{thm:mainthm}.
we can embed the Cayley graph of $\Gamma$ quasiisometrically
into the product of two binary trees.
Since $\H^2$ is quasiisometric to this Cayley graph, we obtain the result.
\end{proof}

\begin{proof} (of Corollary \ref{cor:higherdim})

By a result of Brady and Farb \cite{BF} there exists a quasiisometric embedding of
the hyperbolic space $\H^n$ into the $(n-1)$-fold product of hyperbolic
planes. Actually the proof in \cite{BF}  gives a bilipschitz embedding 
(compare \cite[section 2]{F}).
Combining Corollary \ref{cor:hypplane} with this result, we are done.

\end{proof}

\begin{proof}  (of Corollary \ref{cor:vcohdim})

A recent result of Januszkiewicz and Swiatkowski \cite{JS}
shows the existence of a Gromov-hyperbolic right angled 
Coxeter group $\Gamma_n$
of arbitrary given n, such that the virtual cohomological 
dimension of $\Gamma_n$ is 
$n$. The construction of these groups imply that they have chromatic number
$n$. By Theorem \ref{thm:mainthm} $\Gamma_n$ can be 
quasiisometrically embedded into the product of $n$ binary trees.

\end{proof}

\begin{proof} (of Corollary \ref{cor:hyprank})

We recall the definition of the hyperbolic rank of a metric space.
Given a metric space $M$ consider all locally compact Gromov-hyperbolic
subspaces $Y$ quasiisometrically embedded into $M$. Then
$\rank_h(M) = \sup_Y \dim \partial_{\infty}Y$ is called
the hyperbolic rank. (Compare \cite{BS1} for a discussion of this
notion). 
Let $\Gamma_n$ as in the proof of Corollary \ref{cor:vcohdim} above.
$\Gamma_n$ can be embedded in a
bilipschitz way into the $n$-fold product $T^n$ of the binary tree $T$. Since 
the virtual cohomological dimension of $\Gamma_n$ is $n$, we have
$\dim \partial_{\infty} \Gamma_n = (n-1)$ 
by  \cite{BM}.
Thus $\rank_h(T^n) \geq (n-1)$ by the definition of the hyperbolic rank.
The opposite inequality
$\rank_h(T^n) \leq (n-1)$ follows from standard topological considerations.
 
\end{proof}


\section{Strongly aperiodic tilings of the Davis Complex}

As an application of the methods developed above,
we construct certain aperiodic tilings of
the Davis complex of $X$

We recall the definition of tiling of a metric space from \cite{BW}.
Let $X$ be a metric spcace.
A set of tiles $(\mathcal T,\mathcal F)$ is a finite collection of 
compact $n$-dimensional
complexes $t\in\mathcal T$ and a collection of subcomplexes $f\in\mathcal F$
of dimension $< n$, together with an opposition function $o:\mathcal F\to\mathcal F$,
$o^2=id$.
A space $X$ is tiled by the set $(\mathcal T,\mathcal F)$ if

\begin{itemize}
\item[(1)] $X=\cup_{\lambda}t_{\lambda}$ where each 
$t_{\lambda}$ is isometric to
one of the tiles in $\mathcal T$;
\item[(2)] $t_{\lambda}\setminus\cup_{f\in t_{\lambda}}=Int(t_{\lambda})$ in $X$ for
every $\lambda$;
\item[(3)] If $Int(t_{\lambda}\cup t_{\lambda'})\ne 
Int(t_{\lambda})\cup Int(t_{\lambda'})$ then $t_{\lambda}$ and $t_{\lambda'}$
intersect along $f\in t_{\lambda}$ and $o(f)\in t_{\lambda'}$;
\item[(4)] There are no free faces of $t_{\lambda}$.
\end{itemize}

\begin{Rem}
Strictly speaking the tiling is given by a collection
$\{\phi_{\lambda}\}$ 
of specified
isometries $\phi_{\lambda}:t_{\lambda}\to t$, where $t \in \mathcal T$,
such that for an $n-1$-dimensional face
$\sigma \subset t_{\lambda}$ we have
$\phi_{\lambda}(\sigma) \in \mathcal F$.
Point (3) says the following:
if $\sigma = t_{\lambda} \cap t_{\lambda'}$ is the "common" face,
then $o(\phi_{\lambda}(\sigma)) = \phi_{\lambda'}(\sigma)$
\end{Rem}

Let $g$ be an isometry of $X$ preserving the decomposition structure,
i.e. every
$ t_{\lambda}$ is mapped by $g$ to some
$ t_{\lambda'}$. Then $g$ induces a new tiling
$\{g^*\phi_{\lambda}\}$ with the same tiling set
$(\mathcal T,\mathcal F)$ by
$g^*\phi_{\lambda}=\phi_{\lambda'}\circ g$.
We say that $g$ is an isometry of the tiling, if
$\{g^*\phi_{\lambda}\}=\{\phi_{\lambda}\}$ .
We call a tiling
$\{\psi_{\lambda}\}$ a {\it limit tiling} of
$\{\phi_{\lambda}\}$ , if there exists a sequence $g_i$ of
decomposition preserving isometries of
$X$, such that for every $\lambda$ there exists
$i_0=i_0(\lambda) \in \N$, such that
$\psi_{\lambda}=g_i^*\phi_{\lambda}$ for
all $i \geq i_0(\lambda)$.

\begin{Def}
\begin{itemize}
\item[(1)] A tiling is called {\it  aperiodic}, if the isometry group
of the tiling does not act cocompactly on $X$.
\item[(2)] A tiling is called {\it strongly aperiodic} 
if it has a trivial group of isometries.
\end{itemize}
\end{Def}

\begin{Thm} \label{thm:tiling}
For every finitely generated, right-angled and hyperbolic
Coxeter group $\Ga$ the Davis complex $X$ admits a
strongly aperiodic tiling with finitely many tiles.
In addition every limit tiling is also strongly aperiodic.
\end{Thm}

Since in dimensions 
$2,3,4$ there are right-angled reflection groups
with compact fundamental domain we obtain

\begin{Cor}
The hyperbolic spaces $\H^2,\H^3,\H^4$ admit strongly
aperiodic tilings such that all limit tilings 
are also and strongly aperiodic.
\end{Cor}

Before we prove this result, we recall 
the basic facts concerning the Davis complex. 

\subsection{The Davis Complex}

Let $\Gamma$ be a right angled Coxeter group with generating set $S$ and
let $N=N(\Gamma, S)$ be the nerve of $(\Ga,S)$.
By $N'$ we denote 
the barycentric subdivision of $N$.
The cone $C=\cone N'$  over $N'$ is 
called a {\it chamber}
for $\Gamma$.
The Davis complex \cite{D1} $X=X(\Gamma,S)$ is  
the image of
a simplicial map $q:\Gamma\times C \rightarrow X$ 
defined by 
the following equivalence relation on the vertices:
$a\times 
v_{\sigma}\sim b\times v_{\sigma}$ provided $a^{-1}b\in\Gamma_{\sigma}$.
Here 
$\sigma$ is a simplex in $N$, $\Gamma_{\sigma}$ is the subgroup of
$\Gamma$ 
generated by the vertices of $\sigma$, $v_{\sigma}$ is 
the 
barycenter of $\sigma$.
We identify $C$ with the image
$q(1\times C)$ as a subset of $X$.
The 
group $\Gamma$ acts simplicially on $X$ 
by $\gamma q(\alpha\times x) = q(\gamma\alpha\times x)$
and the orbit space is equal to the chamber $C$.
Thus the Davis complex $X$ is obtained
by 
gluing the chambers $\gamma C$ , $\gamma \in \Gamma$ 
along the boundaries.
Note that $X$ admits an equivariant cell structure with the 
vertices $X^{(0)}$
equal the cone points of the chambers and with 
the 1-skeleton $X^{(1)}$
isomorphic to the Cayley graph of 
$\Gamma$.

The generators $s\in S$ and their conjugates $r=\gamma s \gamma^{-1}$, 
$\gamma\in\Gamma$ are 
called {\it reflections}.
A {\it mirror (or wall)} of a 
reflection $r \in\Gamma$ is the set of 
fixed points $M_r\subset X$ 
of $r$ acting on the Davis complex $X$.
Note that $M_{\gamma s \gamma ^{-1}}= \gamma M_s$.
A wall is indexed by the corresponding
generator $s$. Thus the wall $\ga M_s$ has the
index $s$. If $s \in A$ then we say that the mirror
has {\it $S$-color $a$}. We use the notation $S$-color,
because we will define other colorings of the walls
later. Wall correspond to reflections, thus we
speak also of the $S$-color of a reflection $r$.

We call a boundary wall of a chamber $\ga C$ a face in $X$.
Thus every face in $X$ occurs in exactly two adjacent
chambers. The faces are pieces of walls.

Let $\mathcal M$ be the set of mirrors of the Davis complex,
and let
$\mathcal M_a$ be the mirrors of $S$-color $a$.

\begin{Lem} \label{lem:mirror}
For every generator $s$ in a 
right-angled Coxeter group $\Gamma$ we have
$M_s=\{ q(w \times x) \mid w\in Z_s, x\in St(s,N')\}$,
where $Z_s$ is the centralizer of $s$ in $\Gamma$.
\end{Lem}

\begin{proof} 
$" \ \supset ": \ $ 
Let $w\in Z_s$ and $x\in St(s,N')$, i.e.
$x$ is an affine combination
$x = \sum_{s\in\sigma} x_{\sigma} v_{\sigma}$ one
easily computes
 $$s q(w\times x)=w q(s\times x)= w q(1\times x) = q(w\times x).$$

\noindent $" \ \subset ": \ $ 
Let $z\in M_s$. 
Then $z=q(g\times x)$ for some $g\in\Gamma$ and $x\in cone(N')$. 
The 
condition $s(z)=z$  can be rewritten as $q(sg\times x)=q(g\times x)$.
Hence 
$g^{-1}sg\in\Gamma_{\sigma}$ 
for some simplex $\sigma$ of $N$ and $x\in\cap_{v\in\sigma} St(v,N')$. 
By 
the deletion law 
$s\in\Gamma_{\sigma}$, since the number of $s$ 
in $g^{-1}sg$ is odd.
Hence $s\in\sigma$ and $x\in St(s, N')$.

Let $s_1\dots s_k$ be a reduced presentation of $g^{-1}sg$.
We note that all $s_i\in\sigma$. Since the group is rightangled, $\Gamma_{\sigma}$
is commutative and hence all $s_i$ are different. Let $u_1\dots u_l$ be 
a reduced presentation of $g$. Note that $s_j\ne u_i$ for every $u_i\ne s$,
since $u_i$ appears even number times in the word $u_l\dots u_1su_1\dots u_l$.
Hence $g^{-1}sg=s$, i.e. $g\in Z_s$.
\end{proof}

\begin{Lem} \label{lem:samecolor}

Different 
mirrors of the same $S$-color are disjoint.
\end{Lem}

\begin{proof}
Let 
$a,a'\in A$ and let $M= M_a$ and $M'= M_{a'}$ be the corresponding
mirrors.
Assume that $\gamma_1(M)\cap\gamma_2(M')
\ne\emptyset$.
Therefore $gM\cap M'\ne\emptyset$ where $g=
\gamma_2^{-1}\gamma_1$. 
Let $x\in gM\cap M'$. 
By Lemma \ref{lem:mirror} we have
$x=q(w\times y)=q(gu\times z)$, where
$w\in Z_{a'}$ , $y\in St(t, N')$
and
$u\in Z_a$, $\ z\in St(s, N')$. 
Hence $y=z = \sum_{a,a'\in\sigma} y_{\sigma} v_{\sigma}$.
In particular $a,a'$ are in a common simplex and hence
commute. Thus $a=a'$ since different elements in
$A$ do not commute.
Since $q(w\times y)=q(gu\times z)$, we
have
$w^{-1}gu \in \Gamma_{\sigma}$ for some
simplex $\sigma$ with $a\in \sigma$.
Thus $w^{-1}gu \in Z_a$ which implies
$g \in Z_a$ and $gM_a=M_a$.

\end{proof}

The following Lemma is useful

\begin{Lem} \label{lem:isodavis}
If a group $G$ acts on the Davis complex $X$ of a
Coxeter group $\Ga$ in such a way that $\Ga$ takes
walls to walls and respects the indices of the walls. 
Then $G$ is a subgroup of $\Ga$.
\end{Lem}

\begin{proof}
Let $g \in G$ and aplly $g$ to the base chamber
$C \subset X$. Since $g$ preserves the walls of $X$,
also $gC$ is a chamber, i.e. $gC=\ga C$ as a set.
Since $g$ and $\ga$ are isometries and both
respect the indices of the walls, we clearly
have $g=\ga \in \Ga$.
\end{proof}


\subsection{A tiling by color of the Davis Complex}

Let $X$ be the Davis complex of a right-angled
Coxeter group $\Ga$.

By a {\it coloring} of mirrors by finitely many colors
we mean a function
$\Phi: \mathcal M \to F$, where $F$ is a finite set.
The coloring induces a coloring of the faces of all
chambers, and chambers are glued together respecting the coloring.

\begin{Pro}
Every coloring $\Phi:\mathcal M\to F$ with a finite  $F$ defines a
tiling $(\mathcal T,\mathcal F)$ of the Davis complex $X$  with 
$o(f)=f$ for all $f\in\mathcal F$.
\end{Pro}

\begin{proof}
The set of tiles
$\mathcal T$ is the set of chambers with all posible colorings 
of their faces.
The set of faces $\mathcal F$ is the set of possible
colored faces of the chambers.  Set $o(f)=f$. Then all conditions hold.
\end{proof}

In this situation we can identify
$\mathcal F$ with $F$.
We call such a tiling as {\it tiling by color}.

We note that for any tiling by color 
$\Phi:\mathcal M\to F$ there are limit
tilings.
We can describe them
as follows. The group $\Gamma$ acts (from the right) 
on the compact space
$F^{\mathcal M}$ of
$F$-sequences indexed by $\mathcal M$. 
Every limit point $\Phi'$ of the
sequence
$\Gamma(\Phi)\subset F^{\mathcal M}$
is a limit for some sequence $g_i\in\Gamma$, 
$g_i(\Phi)\to\Phi'$. Then
$\Phi'$
is the limit tiling for $\{g_i\}$.

It is easy to construct a tiling by color of the Davis complex
which is strongly aperiodic. For that it suffices to paint all
walls of the base chamber by different colors and all remaining
wall by another color. Clearly, in this example every limit tiling
is $\Gamma$-periodic. In this section we give a strongly aperiodic
tilings by color of the Davis complex such that all of its limit
tilings are strongly aperiodic.

The proof of 
Theorem \ref{thm:mainthm} gives a coloring of the walls
of the Davis complex in the following way.
Let $M_r$ be a wall, where $r= \ga a\ga^{-1}$ is a reflection for
(say) an element $a\in A\subset S$.
We assume that
$\ell(\ga a) = \ell(\ga)+1$. (If not we replace
$\ga$ by $\ga a$).
Now we define
the color of $\Phi(M_r)$ as the last entry in the sequence
$\psi_a(\ga a)$. 
We have to show that this is well defined:
if we can also write $r = \beta a \beta^{-1}$ with some
$\beta$ and $\ell(\beta) a = \ell(\beta)+1$, then
$\beta^{-1}\ga$ commutes with $a$ and hence does
not contain (in a reduced representation) a letter
of $A$. Note that
$\ga a= \beta a \beta^{-1}\ga$. 
By the definition of
$\psi_a$ we see
that
$\psi_a(\ga a) =\psi_a(\beta a \beta^{-1}\ga) = \psi_a(\beta a)$,
where the last inequality holds since a reduced word 
representing $\beta^{-1}\ga$ does not
contain letters with color $a$.

\smallskip

We investigate what information this coloring gives.
We call two faces of the tiling adjacent, if
the two faces are different but are contained in a
common simplex.

\begin{Def}
A {\it radial gallery of faces} is a sequence 
$F_0,\ldots ,F_k$ of faces, such that there is
a geodesic $\ga_0,\ldots ,\ga_k$ in $\Ga$ starting
from $\ga_0=1$, such that
$F_i$ is a face of the chamber $\ga_i C$.
\end{Def}

\begin{Lem} \label{lem:gallery}
 Let $\Phi$ be the above tiling, and let
$F$ be a face. Then we can reconstruct from the 
properties of the coloring $\Phi$ a radial gallery
$F_0,\ldots ,F_k$ of faces with $F=F_k$.
\end{Lem}

\begin{proof}
Let $F$ be the face between the chambers
$\ga C$ and $\ga a C$ for some generator
$a \in A \subset S$, such that
$\ell(\ga a) = \ell(\ga) +1$.
Now the coloring gives as information the last
entry of
$\psi_a(\ga a)$. 
This last entry contains in particular the diary entry, which
itself consits of the last $\kappa$ entries 
(in the alphabet $\cal{A}$) of the word $W$ which
is the canonical $a$-representation of $\ga$.
The first of this $\kappa$ entries is a (decorated) simplex
$\Delta \subset N(\Ga,S)$.
Note that by definition $s \in \Delta$ if and only if
$\ell(\ga s) < \ell(\ga)$. Choose some $s \in \Delta$ and
let $F'$ be the face between $\ga s$ and $\ga$.
Inductively we will construct a gallery
$F,F',F'',\ldots$ which finally will stop at the base
chamber. Note that the diary entries of all faces of
the base chamber $C$ are $\kappa$ times the symbol $\emptyset$ and a
diary entry is of the form
$(\emptyset,\ldots ,\emptyset)$ if and only if this
entry comes from a face of the basechamber. Thus
the color tells us, when to stop.
In this way we obtain the desired gallery.
\end{proof}

\begin{Rem} \label{rem:diffwalls}
\begin{itemize}
\item[(1)] One easily sees from the proof, that one
can reconstruct from the coloring all possible
geodesics from $\ga C$ to $C$.
\item[(2)] Two different faces of the gallery belong to
different walls.
\end{itemize}
\end{Rem}

Now it is almost immediate that the tiling $\Phi$ is strongly
aperiodic.
Note that as mentioned in the proof of Lemma \ref{lem:gallery}
the faces of the base chamber are characterized by coloring
$\Phi$.
Hence every isometry $g$ of the tiling has to fix the base chamber $C$.
Since $g \in \Ga$ by Lemma \ref{lem:isodavis}, we see $g=1$.

We now consider the case of limit tilings.
Indeed we need some additional properties for the coloring.
We have to make the coloring in a way that "nearby" walls
of the same $S$-color have different colors in $\Phi$.

To make this precise we define

\begin{Def}
Let $M_{r_1}$ and $M_{r_2}$ be two different mirrors with the same $S$-color
(say) $a$.
Then the distance is defined to be the number
$$d(M_{r_1},M_{r_2}) = \inf \{ d(\gamma_1,\gamma_2) + 1 \mid \exists
\gamma_1 ,\ga_2 \in \Gamma, \  a_1,a_2 \in A \ \mbox{with} \ r_i = \gamma_i a_i \gamma_i^{-1}\} $$

\end{Def}

\begin{Rem}

We can view the mirrors $M_{r_i}$ as subsets of the Cayley graph such
that $M_{r_i}$ is the set of midpoints of the edges with
endpoints $\gamma_i $ and $\gamma_i s_i$ where
$\gamma_i s_i \gamma_i^{-1}$ is a representation of $r_i$.
Since by Lemma 
\ref{lem:samecolor} 
two different mirrors with the same color do not
intersect, the distance defined above is exactly the distance
of the mirrors considered as subsets of the Cayley graph.

\end{Rem}

The following result is essential for our construction

\begin{Pro} \label{pro:finite}

For every $D \ge 0$
there exists a map
$f_D:\mathcal M \to G$ where $G$ is a finite
set such that for two reflections
$r_j $, $j=1,2$ which are conjugates of the
same element $s \in S$ the following holds:
$f_D(M_{r_1}) = f_D(M_{r_2})$ implies
$M_{r_1}=M_{r_2}$ 
or $d(M_{r_1},M_{r_2}) \geq D$.

\end{Pro}

\begin{proof}
Let 
$\alpha_1,\ldots ,\alpha_k \in \Gamma$ be the
set of nontrivial elements 
with
$\ell (\alpha_i) \leq 2D +2 $.
Since a Coxeter group $\Gamma$ is
residually finite, there exists a finite
group
$G$ and a homomorphism
$\sigma : \Gamma \to G$
such that
$\sigma (\alpha_i) \neq 1$ for all
$i=1,\ldots,k$ .
For a reflection $r$ we define 
$f_D(M_r) = \sigma(r)$. 
Let $r_1,r_2 $ be as in the assumption and let
$d(M_{r_1},M_{r_2}) < D$.
Then there exits $s \in S$ and $\gamma_j \in \Gamma$ such that
$r_j = \gamma_j s \gamma_j^{-1}$ and
$d(\gamma_1,\gamma_2) \leq D$.
Let
$\tau = \gamma_1^{-1}\gamma_2$,
hence $\ell(\tau) \leq D$.
By assumption
$r_1^{-1}r_2 \in ker(\sigma)$.      
Since
$r_1^{-1}r_2=\gamma_1s\tau s \tau^{-1}\gamma_1^{-1}$
we have
$s \tau s \tau^{-1} \in ker (\sigma)$.
Since
$\ell(s \tau s \tau^{-1}) \leq 2 D + 2$
we have by construction that
$s \tau s \tau^{-1}$ is trivial and
hence
$\tau$ commutes with $s$ which implies that
$r_1 = r_2$.
\end{proof}

\begin{Pro}
For every finitely generated, right-angled and hyperbolic
Coxeter group $\Ga$ the Davis complex $X$ admits a
tiling $\Psi$ by colors with finitely many tiles such
that $\Psi$ and every limit tiling of 
$\Psi$ is also strongly aperiodic.
\end{Pro}

\begin{proof}
We have to construct a coloring of the walls
$\mathcal M$ by finitely many colors, such that the corresponding
tiling is strongly aperiodic and also all limit tiling
are strongly periodic.
First we associate to every wall
$M_r$ the "color"
$\Phi(M_r)$.
This gives already information
about

\begin{itemize}
\item the $S$-color of $M_r$
\item the diary entries in order that
Lemma \ref{lem:gallery} holds.
\item the Morse Thue decoration.
\end{itemize}

We need in addition a second "decoration", 
which we will call the 
{\it distance decoration},
of the walls with the
following property:
If $M$ and $M'$ are walls with the same
$S$-color and the same distance decoration,
then 
$$M = M' \ \ \mbox{or} \ \ d(M,M') > 2 C_h +3,$$ 
where
$C_h$ is the hyperbolicity constant of
$\Ga$.
It follows directly from Proposition \ref{pro:finite}
that we can find such a decoration with a finite
decorating set.

If we combine the coloring $\Phi$ with this
additional decoration, we obtain a coloring
$\Psi:\mathcal M \to F$, where $F$ is some finite set.
Thus the induced tiling has again only finitely many tiles.
As already shown above the tiling is strongly aperiodic and
we have to consider a limit tiling.
In the case of a limit tiling 
Lemma \ref{lem:gallery} then implies that
given a face $F_0$ of the tiling, we can construct a gallery
$F_0,F_{-1},F_{-2} \ldots$ starting from $F_0$.
We can assume that this gallery does not end, since
we are in a limit case.
Note that this
gallery can be considered as a geodesic 
$\nu:\{0,1,2,\ldots \} \to C(\Ga,S)$,
where every point
$\nu(i)$ correpond to a face of the tiling
or an edge in $C(\Ga,S)$.
Assume that there is an isometry
$g$ of the tiling.
Then $g \nu$ is also a geodesic in $C(\Ga,S)$.
Since
$\nu$ and $g \nu$ are limits of isometric images of
geodesics segments starting
at the base chamber, the geodesics
$\nu$ and $g \nu$ are asymptotic.
Thus there exists a "shift constant" $k$ (depending on
$g$) 
such that
$d(\nu(i),g \nu (i+k)) \leq C_h$ for all $i$ large enough, where 
$C_h$ is the hyperbolicity constant of $\Ga$.
By shifting the initial point of the geodesic and maybe interchanging
the role of $\nu$ and $g \nu$, we can
assume w.l.o.g. that $k\geq 0$ and
$d(\nu(i),g \nu(i+k)) \leq C_h$ for all $i \ge 0$.
There exists some $S$-color (say) $a$, such that this color
occurs infinitely often in this sequence. 
Now fix in addition a possible decoration.
Let
$j_0, j_1, \ldots $ be the subsequence such that the
$\nu (j_i)$ are the faces of $S$-color $a$ and this given
decoration. Note that 
this sequence may be finite.
By Remark \ref{rem:diffwalls} (2) the walls belonging to
the faces $\nu (j_i)$ are different.
By the properties of the distance decoration
there is a sequence $k_i$, with
$k_i \le j_{i} \le k_{i+1}$, such that
the faces $\nu (k_i)$ and $\nu(k_{i+1})$  have distance $\ge C_h$ from
the walls belonging to the face
$\nu (j_i)$ and these points lie on different sides of this wall.
By triangle inequality also the points
$g \nu(k_i)$ and $g \nu (k_{i+1})$ are on different sides
of this wall which implies that $g \nu$ has to cross the wall.
Thus the geodesic $g \nu$ also crosses all the faces with $S$-color
$a$ and the given decoration which occur on the geodesic $\nu$.
This argument holds for all possible decorations. Hence
all walls of $S$-color $a$, which are intersected by $\nu$ are
also intersected by $g \nu$ and 
(since walls of the same $S$-color do not intersect) in the same
order.
Thus the sequence of these faces on $\nu$ is a shift these faces
on $g \nu$. By the properties of the Morse Thue decoration, there are
no propper shifts. This implies that
$g$ leaves these walls invariant. This implies $k=0$ and
hence the displacement of $g$ is $\leq C_h$ on $\nu (i)$.
Thus there is a chamber of $X$ which is moved by $g$ at most
by distance $C_h+1$. Since $g$ preserves also the distance decoration,
all the walls of this chamber are invariant under $g$ which implies
that $g$ fixes the chamber and hence (since the coloring also
contains the information of the element $s \in S$ of
the wall $M_r$, where $r$ is a conjugate of $s$),
$g$ is the identity. 

\end{proof}


\subsection{Balanced Tilings} \label{subsec:balanced}

Let $(\mathcal T,\mathcal F)$ be a set of tiles. 
A function $w:\mathcal F\to\bf Z$ is
called {\it a weight function} if
$w(o(f))=-w(f)$ for every $f\in\mathcal F$.
We recall a definition from [BW].
\begin{Def}\label{def:Definition 6.16}
A finite set of tiles $(\mathcal T,\mathcal F)$ 
is unbalanced if there is a weight
function $w$ such that
$\sum_{f\in t}w(f)>0$ for all $t\in\mathcal T$.

It is called semibalanced if $\sum_{f\in t}w(f)\ge 0$ for all 
$t\in\mathcal T$.
\end{Def}

We call a set of tiles 
{\it strictly balanced} if for every nontrivial
weight function $w$ there
are tiles $t_+$ and $t_-$
such that $\sum_{f\in t_+}w(f)>0$ and 
$\sum_{f\in t_-}w(f)<0$. A tiling
called
{\it unbalanced (strictly balanced)}
if the corresponding set of tiles is unbalanced (strictly balanced).

In [BW] aperiodic tilings of some nonamenable metric spaces (such as the
Davis complexes
of hyperbolic Coxeter
groups) are constructed where the aperiodicity follows from the fact that
they are unbalanced.
Here we show,
how one can modify a tiling $\Phi$ by color of the Davis
complex in a way, that the new tiling is strictly balanced and also all
limit tilings are strictly balanced.
We start from a tiling by color, i.e. a function
$\Phi:\mathcal M \to F$ where $F$ is a finite set.

We now associate in addition to every wall an orientation.
A wall divides the Davis complex into two components. Roughly speaking the
orientation says which of the components is left and which is right.
The orientation of the walls defines a new tiling
$(\mathcal T',\mathcal F')$ where the new faces have in addition
a sign $+$ or $-$.
Thus $\mathcal F' = \mathcal F_{+} \cup \mathcal F_{-}$, where
$\mathcal F_{+}$ and  $\mathcal F_{-}$ are copies of $\mathcal F$.
The face $f \in t_{\lambda}$ has sign $+$, if
$Int(t_{\lambda})$ is left of the wall and sign $-$, if
$Int(t_{\lambda})$ is right of the wall.
In this case the opposition function
$o:\mathcal F'\to \mathcal F'$ maps
$f_+$ to $f_-$.
Geometrically this means that we deform all faces of given color and given sign
in the same direction by the same pattern. For the faces of the same color but
of opposit sign we take the opposite deformation. We call such tiling
a {\it geometric resolution} of a tiling by color.
This new tiling is not any more a tiling by color.

Note that in the Davis complex a wall has a canonical orientation,
by deciding that the base chamber is in the left component.
Thus we can indicate the choosen orientation itself by a sign.
A wall gets the sign $+$, if the orientation of the wall is the
canonical one and $-$ otherwise.

The following is obvious.
\begin{Lem}
(1)Assume that a coloring $\Phi:\mathcal M\to F$ 
defines a strongly aperiodic tiling of $X$. Then
any geometric resolution is also strongly
aperiodic.

(2) If every limit tiling of $\Phi$ is strongly aperiodic,
then every limit tiling of any
geometric resolution is strongly
aperiodic.

\end{Lem}

We formulate and prove the main result of this section.

\begin{Thm} \label{thm:balanced}
For every coloring $\Phi:\mathcal M\to F$ with the
property that walls of the same color do not intersect, there is a
strictly balanced geometric resolution of this tiling. In addition also
every limit tiling of this resolution is strictly balanced.
\end{Thm}

\begin{proof}
In the first step of the proof we
construct a strictly balanced geometric resolution of
$\Phi$.
Consider the set of mirrors  $\mathcal M_f=\Phi^{-1}(f)$ 
of the same color $f$.
Since walls of the same color do not intersect
by assumption, they are ordered by level from the base chamber.
(The level is defined by induction.
If one removes the mirrors $\mathcal M_f$ from $X$, the space is
devided into components.
Mirrors that bound the
component of the base chamber are of level one. Then drop mirrors 
of level one
and repeat the procedure to get new mirrors of level one and call 
them of level
two and so on). 
We give the mirrors signs in alternate fasion by the level:
$+-+-+-\dots$ . Thus, elements of $\mathcal M_f$
have in alternate ways the $+$-orientation and the $-$-orientation.

The choosen orientation induces a geometric resolution of the tiling.
We show that this resolution is strictly balanced.

Let $w:F_+\cup F_-\to\Z$ be a nontrivial 
weight function with $w(f_+)=-w(f_-)$. 
We show that there is a chamber $C_+$ such that
$$
\sum_{M\in C_+}w(M)>0
$$
and
there is a chamber $C_-$ such that
$$
\sum_{M\in C_-}w(M)<0.
$$

We show the first.

Since the function
$w$ is nontrivial there exists a face
$f^0$ which is the common face of two
adjacent complexes
$t_{\lambda}$ and $t_{\lambda'}$ such that
$w(f^0)\neq 0$.
(By a slight abuse of notation we here identify 
the "real face" $\sigma = t_{\lambda}\cap t_{\lambda'}\subset X$ 
with the color
$\Phi(\phi_{\lambda}(\sigma))=\Phi(\phi_{\lambda'}(\sigma))\in \mathcal F$.)

Now there are four cases corresponding to the
parity of the sign of
$w(f^0_+)$ and the orientation of the
mirror $M_{f^0}$ defined by the face $f^0$.
All this cases are treated in a similar way.
We discuss only one, and to be fair not the easiest of these
cases:
let $w(f^0_+) > 0$ and the orientation of the wall
$M_{f^0}$ be negative.

By the convention how we use the words "left" and
"right" this means that the base chamber is
right of the wall $M_{f^0}$.
Now
$f^0= t_{\lambda'}\cap t_{\lambda}$ where we assume that
$t_{\lambda}$ is right, and $t_{\lambda'}$ is left of the
wall. Thus
$f^0$ considered as a face of
$t_{\lambda}$ has sign $-$ and considered as a face of
$t_{\lambda'}$ has sign $+$.
Now we choose some number
$k \in \N$ such that every geodesic
from
$Int(t_{\lambda})$ to $Int(C)$ intersects $\leq k$ mirrors.
(i.e. $k$ is larger than the combinatorial distance
of $t_{\lambda'}$ and the base chamber.)

Now consider any color $f$.
Then there are three possibilities:

(a) $w(f_+) > 0$ and $w(f_-)< 0$; then  we call
$f$ an {\it even} color for $w$,

(b) $w(f_+) < 0$ and $w(f_-)>0$; then  we call $f$ 
an {\it odd} color for $w$.

(c) $w(f_+) = w(f_-) =0$; then  we call $f$ 
a $0$-color for $w$.

\noindent Let 
$$\mathcal M_{ev}^{2k} = \{M_f | f \ne f^0,\  f \ \ \mbox{even for} \ w\  \mbox{and level of}\  
f = 2k\}$$
$$\mathcal M_{odd}^{2k+1} = \{M_f | f \ne f^0,\  f \ \ \mbox{odd for} \ w\  \mbox{and level
of}\ f =2k+1\}$$
$$\mathcal M_0^{k+1} = \{M_f | f \ \ \mbox{0-color for} \ w\  \mbox{and level
of}\ f =k+1\}$$

Claim 1:  {\it The set of mirrors 
$\mathcal M_{ev}^{2k}\cup \mathcal M_{odd}^{2k+1} \cup \mathcal M_0^k
\cup \{ M_{f^0}\}$
bound a bounded set $D$ containing the chamber $t_{\lambda'}$.} 
 
Clearly this set of mirrors bounds a convex set in the
Hadamard space $X$. If the component containing
$t_{\lambda'}$ is unbounded, then there is a ray from
$t_{\lambda'}$
to the visual boundary which does not intersect any of our mirrors.
Since we have only finitly many colors, there is a color
$f$ such that this ray intersects infinitely many
mirrors of this color. By the choice of $k$ the first of these intersected
mirrors has level $\le k+1$ and mirrors of the same color
and level are intersected at most twice by convexity.
Thus one of these intersected mirrors is contained in
our set of mirrors.
Contradiction.

Claim 2: {\it If $f$ occurs as face of a tile 
$t_{\mu} \subset D$ such that
$f \subset \partial D$, then $w(f) \geq 0$.} 

To prove Claim 2 we consider the cases

\noindent (i) If $f =f^0$, then $f$ as a face of $t_{\mu}$ has
the orientation $+$ since $t_{\mu}$ lies on the same side of
$M_{f^0}$ as $t_{\lambda'}$. Thus
$w(f)=w(f^0_+)>0$.

\noindent (ii) If $f$ is a $0$-face , then anyway
$w(f) = 0$.

\noindent (iii) Let $f$ be an even color for $w$. 
Then $f$ is contained in a wall of $\mathcal M_{ev}^{2k}$.
By the choice of $k$ the tile
$t_{\mu}$ is on the same side of this wall as the base chamber.
The wall has level $2k$ hence orientation $+$, thus the basechamber and hence
$t_{\mu}$ is on the left side of this wall.
This implies that $f$ considered as a face of
$t_{\mu}$ gets the sign $+$. Hence
$w(f) = w(f_+) > 0$ since $f$ is an even color for $w$.

\noindent (iv) A similar argument aplies for $f$, if
$f$ is an odd color for $w$. This proves Claim 2.

\smallskip

According to the Claim 1 
we have $D=\cup_{i=1}^kC_i$ where $C_1,\dots, C_k$ is a finite
collection of chambers. Then 
$$
\sum_{i=1}^k\sum_{f'\in C_i}w(f')=\sum_{f'\in\partial D}w(f')\geq 0.
$$
by Claim 2. Since also
$f^0_+$ is in the last set of faces we see that the expression is
indeed $>0$.
Therefore, $\sum_{f'\in C_i}w(f')>0$ for some $i$.

\smallskip

\noindent To obtain the chamber $C_-$, we make a similar construction.

This finishes the proof of the first step.
Thus we have constructed a strictly balanced geometric resolution
of $\Phi$.

Actually the proof of the first step shows more:
If we choose for any given color $f$ an orientations of
the walls
$\mathcal M_f$ in alternate way $+-+-\ldots\ $ {\it or}
$\ -+-+ \ldots$ (and maybe for different colors in a different way),
then the resulting geometric resolution is strictly balanced.
(This more general result follows from some obvious
modifications of the above proof).
Let us call such a choice of orientations
an {\it allowed} orientation of walls.
The levels of walls depend on the base
chamber. If we define levels with respect to a different chamber,
they are changed. The change has the property, that either
all parities of the levels are preserved or all are changed.
Thus, whether an orientation of the walls is allowed (or not),
does not depend on the choice of the basechamber.
As a consequence we have the following:
if the orientation of a tiling by color
$\{\phi_{\lambda}\}$ is allowed, then
also the orientation of
$\{g^*\phi_{\lambda}\}$ is allowed.
Thus also all limit tilings of the tiling
constructed in step 1 are strictly balanced.

\end{proof}

 In twodimensional jigsaw tiling puzzles a geometric
 resolution is usually realized by adding rounded tabs out on the sides
 of the pieces with a corresponding blank cut into the intervening
 sides to receive the tabs of adjacent pieces.
 This procedure destroyes the convexity of the pieces.
 We show that in the case of the hyperbolic plane
 $\H^2$ we can modify this construction to obtain aperiodic and 
 strictly
 balanced tilings with convex tiles.
 Compare also the papers \cite{MaMo}, \cite{Mo}.
 
 \begin{Thm}
 (1) For every $n\ge 3$ there is a strictly balanced strongly
aperiodic tiling of $\H^2$ by convex $2n$-gons
with finitely many tiles such that every limit tiling is strongly
aperiodic.

(2) For every $n\ge 3$ there is a finite set of tiles
$(\mathcal T,\mathcal F)$ that consists of convex $2n$-gons
with a strongly aperiodic tiling of $\H^2$ with strongly aperiodic
limit tilings such that every $(\mathcal T,\mathcal F)$-tiling
of $\H^2$ is aperiodic.
 
 \end{Thm}
 
 \begin{proof}
 (1)
 Identify $\H^2$ with the Davis complex
 for the right-angled Coxeter group generated
 by the reflections at a regular right angled
 $2n$-gon. This is a 2-colored group. We call
 these $S$-colors $a$ and $b$.
 We fix a strongly aperiodic tiling $\Phi$ by color as constructed
 above and an orientation of the walls in order
 that the geometric resolution of the tiling is strictly
 balanced. 
 We now define a modification of the tiling defined by
 $\Phi$. Consider a vertex of the Davis complex. This is a point
 where some $a$-wall intersects some $b$-wall. Thus due to the 
 map $\Phi$  the point correspond
 to a unique pair of colors
 $(a_i,b_j)\in F\times F$, where $F$ is the finite image of
 $\Phi$ and $a_i$ (resp. $b_j$) indicates that the corresponding
 $S$-color is $a$ (resp. $b$).
 The additional information about the orientation of the walls
 give a well defined lokal coordinate system around the vertex
 (by deciding that the positive quadrant is the quadrant which is
 right to both walls, and by deciding that the $a$-wall correspond
 to the first and the $b$-wall to the second coordinate.)
 We now move the vertex by a small amount using these local
 coordinates. This move should only depend
 on $(a_i,b_j)$ and not on the vertex. 
 E.g. we can choose a small $\mu(a_i,b_j) > 0$ and
 move the point to a distance $\mu(a_i,b_j)$ into the 
 direction of the diagonal of the positive quadrant.
 After this deformation we obtain a finite number of new convex tiles, which
 (for generic deformations) only allow tilings of $\H^2$
 compatible with the matching rules defined by $\Phi$ and the
 orientation of the walls. Roughly speaking the matching
 rules are now encoded in the length of the sides and the
 angles of the tiles. Our original aperiodic tiling
 is deformed to an aperiodic tiling with the desired properties.

  (2) We start with the tiling by color $\Phi:\mathcal M\to F$ 
defined earlier and take an
unbalanced geometric resolution.
One of the way to do it is to assign + for passing through a wall
from the component that contains the base chamber. Then for every
chamber $C$ the faces whose walls do not separate $C$ and the base
chamber obtain the sign +, all other -. In view of hyperbolicity
the number of + is greater than - for every $C$. We define a
weight function by sending a positive face to $+1$ and a negative
face to $-1$. This geometric realization has all required
properties except the last one by (1). Since the set of tiles
$(\mathcal T,\mathcal F)$ is unbalanced, in view of the following
proposition any other tiling by $(\mathcal T,\mathcal F)$ is
aperiodic.
\end{proof}

\begin{Pro} \label{pro:Proposition 1.19}
Let $(\mathcal T,\mathcal F)$ denote the set of 
tiles of a geometric realization
of a tiling by color of the Davis complex $X$
of a Coxeter group $\Gamma$ supplied with a left 
$\Gamma$-invariant
metric. Suppose that the set of tiles
$(\mathcal T,\mathcal F)$ is unbalanced. 
Then any $(\mathcal T,\mathcal F)$-tiling of $X$
is aperiodic.
\end{Pro}

\begin{proof}
This proposition can be derived formally from Proposition 4.1 [BW].
Since the proof there has some omissions we
present a proof below.

Let $G$ be a group of isometries of a 
$(\mathcal T,\mathcal F)$-tiling of
$X$ which acts cocompactly. Since the group of isometries of $X$
is a matrix group, by Selberg Lemma $G$ contains a torsion free
subgroup $G'$ of finite index. Then the orbit space $X/G'$ is
compact and admits a $(\mathcal T,\mathcal F)$-tiling 
(Note that by taking
$X/G$ as in [BW] we cannot always obtain a tiling because of free
faces). Then we obtain a contradiction:
$$
0<\sum_{t\in G/G'}\sum_{f\in t}w(f)=\sum_{f\in X/G'}w(f)+w(o(f))=0.
$$
\end{proof}


\bigskip
\begin{tabbing}

Alexander Dranishnikov,\hskip10em\relax \= Viktor Schroeder,\\ 

Dep. of Mathematics,\>
Institut f\"ur Mathematik, \\

University of Florida,\> Universit\"at Z\"urich,\\
444 Little Hall, \>
 Winterthurer Strasse 190, \\

Gainesville, FL 32611-8105\>  CH-8057 Z\"urich, Switzerland\\

{\tt dranish@math.ufl.edu}\> {\tt vschroed@math.unizh.ch}\\

\end{tabbing}
\end{document}